\documentclass[a4paper]{article}
\pdfoutput=1
\usepackage[latin1]{inputenc}
\title{On convex regression estimators}
\newcommand{\CONICET}{%
   Consejo Nacional de Investigaciones Científicas y Técnicas}
\newcommand{\UNL}{Universidad Nacional del Litoral}
\author{
   Néstor Aguilera\thanks{\CONICET, and \UNL, Argentina}
   \and
   Liliana Forzani\footnotemark[1]
   \and
   Pedro Morin\footnotemark[1]
   }

\usepackage{ifthen}
\newboolean{hyper}\setboolean{hyper}{true} 
\usepackage{cmap}       
\usepackage[T1]{fontenc}
\usepackage{lmodern}
\usepackage{textcomp}   
\usepackage{microtype}
\usepackage{natbib} 
\usepackage{amsmath,amsthm}
\usepackage{amssymb,mathrsfs}
\usepackage{enumerate}
\usepackage{graphicx}
\usepackage[small]{subfigure}
\usepackage{color}
\usepackage{xspace}
\usepackage{setspace}
\ifthenelse{\boolean{hyper}}{%
   \usepackage[%
      pdftex=true, 
      colorlinks,
      linkcolor=blue,
      anchorcolor=blue,
      citecolor=blue,
      urlcolor=blue
      ]{hyperref}%
   \usepackage{hypcap}%
   }{
   \usepackage{url}%
   }
\ifthenelse{\isundefined{\backref}}{}{
   \renewcommand*{\backref}[1]{} 
   \renewcommand*{\backrefalt}[4]{%
      \ifcase #1 {\textcolor{red}{sin citas}}   
      \or (\textit{pág.~#2})        
      \else (\textit{págs.~#2})     
      \fi}
   }
\usepackage{verbatim}

\theoremstyle{plain}
   \newtheorem{thm}{Theorem}
   \newtheorem{lem}[thm]{Lemma}
   \newtheorem{cor}[thm]{Corollary}
   \newtheorem{proc}[thm]{Procedure}
\theoremstyle{remark}
   \newtheorem*{rem}{Remark}
\newcommand{\NN}{\mathbb{N}}
\newcommand{\RR}{\mathbb{R}}
\newcommand{\convfs}{\mathscr{C}} 
\newcommand{\conve}[1][{}]{f^c_{#1}}
\newcommand{\convn}{\conve[n]}
\newcommand{\convl}{\mathscr{L}} 
\DeclareMathOperator{\PP}{\mathbb{P}}
\newcommand{\prob}[1]{\PP(#1)}
\newcommand{\Prob}[1]{\PP\!\left(#1\right)}
\DeclareMathOperator{\EE}{\mathbb{E}}
\newcommand{\esp}[1]{\EE(#1)}
\newcommand{\Esp}[1]{\EE\!\left(#1\right)}
\newcommand{\as}{a.s.\xspace}
\newcommand{\norm}[1]{\Vert#1\Vert}

\newcommand{\abs}[1]{\vert#1\vert}

\newcommand{\card}[1]{\#(#1)} 
\newcommand{\Lip}{\text{\upshape Lip}}
\newcommand{\cube}{Q}
\newcommand{\grid}{\mathcal{M}}
\makeatletter
\newcommand{\ifsuba}[2]{#1_{n\@ifnotempty{#2}{,#2}}} 
\newcommand{\ifsupa}[2]{#1^n\@ifnotempty{#2}{_{#2}}} 
\newcommand{\ifsubb}[2]{#1\@ifnotempty{#2}{_{#2}}}
\makeatother
\newcommand{\hn}[1][{}]{\ifsuba{h}{#1}} 

\newcommand{\xx}[1][{}]{\ifsubb{X}{#1}}
\newcommand{\yy}[1][{}]{\ifsubb{Y}{#1}}
\newcommand{\ee}[1][{}]{\ifsubb{e}{#1}}
\newcommand{\dd}{d} 

\newcounter{savectr}
\newcommand{\savectr}[1]{\setcounter{savectr}{\value{#1}}}
\newcommand{\backctr}[1]{\setcounter{#1}{\value{savectr}}}

\newcommand{\tiponame}{T-}
\newenvironment{tipo}{%
   \begin{enumerate}[{\bfseries{\tiponame}1.}]}{%
   \end{enumerate}}
\newcommand{\condname}{H-}
\newenvironment{cond}{%
   \begin{enumerate}[{\bfseries{\condname}1.}]}{%
   \end{enumerate}}
\ifthenelse{\boolean{hyper}}{
   \newcommand{\tiporef}[1]{\hyperref[#1]{{\bfseries{\tiponame\ref*{#1}}}}}%
   \newcommand{\condref}[1]{\hyperref[#1]{\bfseries{\condname\ref*{#1}}}}%
   }{
   \newcommand{\tiporef}[1]{\bfseries{\tiponame\ref{#1}}}%
   \newcommand{\condref}[1]{\bfseries{\condname\ref{#1}}}%
   }
\ifthenelse{\boolean{hyper}}{
   \newcommand{\tref}[2]{\hyperref[#2]{#1~\ref*{#2}}}%
   \newcommand{\eref}[2]{\hyperref[#2]{#1~\eqref{#2}}}%
   }{
   \newcommand{\hyperref}[2][]{#2}
   \newcommand{\tref}[2]{#1~\ref{#2}}%
   \newcommand{\eref}[2]{#1~\eqref{#2}}%
   }

\newcommand{\qand}{\quad\text{and}\quad}
\newcommand{\qhull}{\textsf{QHULL}\xspace}
\newcommand{\matlab}{\textsf{MATLAB}\xspace}
\newcommand{\rutina}[1]{\textsf{#1}\xspace}

\begin{document}
\bibliographystyle{sjs-modif}
\maketitle
\thispagestyle{empty}

\begin{abstract}
A new nonparametric estimator of a convex regression function in any dimension is proposed and its convergence properties are studied.
We start by using any estimator of the regression function and we \emph{convexify} it by taking the convex envelope of a sample of the approximation obtained.
We prove that the uniform rate of convergence of the estimator is maintained after the convexification is applied.
The finite sample properties of the new estimator are investigated by means of a simulation study and the application of the new method is demonstrated in examples.
\end{abstract}


\section{Introduction}\label{sec:intro}

In the nonparametric regression model
\begin{equation}\label{equ:model}
   \yy[n] = f(\xx[n]) + \ee[n], \quad n = 1, 2,\dotsc,
\end{equation}
where $\yy[n]\in\RR$, $\xx[n]\in\RR^d$ and $\ee[n]$ is an error term, it is not uncommon to have strong presumptions on properties of $f$---such as monotonicity, convexity or concavity---which should be taken into account.

Typical examples appear in economics (indirect utility, production or cost functions), medicine (dosage-response experiments) and biology
(growth curves).

A much studied case is the instance of a monotone regression function for $d = 1$, estimated by using least squares
(see, e.g.,
\citealp{Br55,Mu88}, 
and
\citealp{BBBB72} 
or \citealp{RoWrDy88} 
for a summary of this work).
For convex (concave) regression
\citet{Hi54} 
proposed to use convex least square estimates, and
\citet{HaPl76} 
proved their consistency.
Algorithms for computing these estimates were developed by
\citet{Wu82} 
and
\citet{FrMa89}, 
and the rate of convergence was derived by
\citet{Ma91}. 
Later
\citet{GJW01} 
derived the asymptotic distribution of the estimator at a fixed point of positive curvature.
In all of these works the estimates hold pointwise.

Still in one dimension, one can avoid the complications of least squares techniques and use more conventional smoothing methods when $f$ is convex (or concave), as shown by
\citet{BD07}. 
Using the fact that a differentiable function is convex (concave) if the derivative is increasing (decreasing), they propose to first smooth the data using any constrained nonparametric estimate (kernel type, local polynomial, series or spline estimator), then compute the derivative of the smooth function thus obtained, which is isotonized and finally integrated to recover a convex estimation.
As mentioned above, the isotonization of a function is something that has already been mastered in the non-parametric literature, and using those results the rates of convergence obtained by them are the usual in non-parametric regression.

Unfortunately this technique can only be used in one dimension and with smooth convex functions and cannot be extended to higher dimensions,
since there is no such simple characterization of convexity in $\RR^d$ for $d > 1$.

As far as we know, little has been done in higher dimensions.
\citet{SHH05} (see also~\citealp{HSHKH05}) present a multivariate data smoothing method using a linear program (for the $\ell^1$ and $\ell^\infty$ norms) or quadratic program (for the $\ell^2$ norm).
\citet{SCK06} develop an approximation method based on multivariate adaptive regression splines (MARS).
But none of these articles present convergence results.

We propose here a simple and fast method that can be used in any dimension and applied to any convex function, even if not too smooth.
Like Birke and Dette, we start by using any approximating scheme on the data, but then we use a \emph{convexification} step, consisting in taking the convex envelope of the approximating function just obtained.
This last step can be done very quickly by current software such as \qhull \citep{qhull}, and the uniform rate of convergence of the approximation technique is maintained after the convexification is applied.

More precisely, we obtain uniform error estimates, and the rate of convergence of the convex estimator is the same as that of the original estimator, thereby showing that the convexification step adds basically no further errors to the estimating step.

The paper is organized as follows.
In \tref{Section}{sec:review} we briefly review fundamental smoothing techniques.
In \tref{Section}{sec:convexification} we show theoretical results on the
convexification
step, and how the error estimates for the convex estimate are derived from the smoothing step.
Finally, in \tref{Section}{sec:numerical} we apply these techniques to approximate several problems in dimensions $d = 1$ and $d = 2$.

\section{The smoothing step: review of the literature}
\label{sec:review}

As we have already pointed out, our method of convexification inherits the $L^\infty$ rate of convergence from whichever smoothing process is chosen for the \eref{model}{equ:model}.
We think it is appropriate, then, to briefly review rates of convergence in $L^\infty$-norm for some of the possible choices for such a process when no monotonicity or convexity assumptions are made on $f$.

Most of the approximation techniques with known rates of convergence are of the so called \emph{smoothing} type, where a variable kernel is used, and we will focus our attention on these.

It should be noted that since there are many different schools and people involved, here we can give only partial references, leaving out several meaningful results available in the literature.

Perhaps the first ones to consider these problems were \citet{De78} and \citet{MR515689}.
Devroye considered the Nadaraya-Watson regression estimator and proved the uniform convergence (without rates) for independent data, with fixed or random predictors belonging to $\RR^d $, whereas Schuster and Yakowitz considered more general kernels in one dimension, establishing orders of convergence in probability.
Later these results were extended by several authors,
among them \citet{MR721221} and \citet{MR920343}.
They extended the result to non-independent data and Collomb was the first to give strong rates for uniform convergence.
Further results on uniform convergence rates for different settings such as robust estimation and other kind of non-independent data were given by \citet{MR866288}, \citet{Roussas}, \citet{Boente-Fraiman}, \citet{MR1150335} and \citet{Tran}.
Extensions to spline estimators were given by \citet{MR2387772}, and to uniform choice of bandwidth by \citet{MR2195639,MR1744994}, \citet{Dony:2008}, \citet{MR2261055}, \citet{MR2543588}, and \citet{DEM:2006} (see also the references therein).

The asymptotic distribution of the maximal deviation between a non-par\-a\-metric regression estimator and the true regression was first considered by \citet{Jh82}, extending to the regression context the results by \citet{MR0348906} and \citet{MR0428580} on density estimation.
For the case $d = 1$ and random predictors, \citeauthor{Jh82} showed---under some regularity assumptions---the $L^{\infty}$ asymptotic distribution of the kernel regression estimator, which allowed him to give uniform confidence intervals for the regression estimator.
This result was extended by \citet{MR768499} to other kernel estimators and by \citet{Ha89} to general estimators defined implicitly, as for example $M$-smoothers and local polynomial estimators.
As far as we know these results were not extended to higher dimensions or non-independent data.

\section{A convex estimator and its convergence}
\label{sec:convexification}

Let us assume that the variables $\xx[n]$ in the \eref{model}{equ:model} take values on a bounded closed convex set $\cube\subset\RR^d$, and that
$f\in\convfs$, where $\convfs$ is the set of (finite real valued) convex functions defined on $\cube$.

$\cube$ need not be polyhedral, but assuming its boundary is smooth except for a finite set of ``corners'', in practice we may approximate it by a polyhedron.
Thus, from now on, for simplicity we will assume that $\cube$ is a polyhedron, and therefore it is the convex hull of its finite set of vertices.
In particular, we assume that $\cube$ is compact.

Let us assume that $f_n$ is an estimator of $f$, defined in all of $\cube$.
To fix ideas, we may think that $f_n$ is obtained by considering the points $(\xx[i],\yy[i])$, $i = 1,\dots,n$, by some procedure such as smoothing.
Our purpose is to derive from $f_n$ another estimator which is also convex.

To do so, we consider a finite set $\grid_n\subset\cube$ such that the convex hull of $\grid_n$ is $\cube$.
The number of points in $\grid_n$ need not be $n$ and the points in $\grid_n$ might be completely unrelated to $\{\xx[i] : i\in\NN\}$.

We now let $\convl_n$ be the set of ``convex functions below $f_n$ on $\grid_n$'',
\[
   \convl_n = \bigl\{\psi\in\convfs: \psi(x) \le f_n(x) \text{ for all $x\in\grid_n$}\bigr\},
\]
and define the \emph{convex estimator} $\convn$, associated with the estimator $f_n$ and the set $\grid_n$ by
\begin{equation}\label{equ:convn:def}
   \convn = \sup \,\{\psi : \psi\in\convl_n\}.
\end{equation}

Since $\grid_n$ contains all the vertices of $\cube$, it is easy to see that $\convn$ is well defined on $\cube$ and that $\convn\in\convfs$.
Furthermore, $\convn$ is piecewise linear, determined by the maximum of hyperplanes.
In particular:

\begin{lem}\label{lem:convn:Ln}
$\convn\in\convl_n$.
\end{lem}

As $\convn$ is the ``lower part'' of the convex hull of the set $\{(x,f_n(x)) : x\in\grid_n\}$, we may take advantage of any of a number of algorithms for finding convex hulls in $\RR^d$.
For instance, \qhull \citep{qhull} finds the convex hull of a finite set of points in any number of dimensions, and is really fast for dimensions $d \le 4$.

We are led to the following procedure for constructing a convex estimator $\convn$ of $f$:

\begin{proc}\label{proc:convn}
Given $\xx[i]$ and $\yy[i]$ ($i = 1,2,\dotsc$):

\begin{enumerate}[Step 1.]
\item
(Smoothing)
Construct an estimator $f_n$ of $f$, for instance through a smoothing procedure using the values $\xx[i]$ and $\yy[i]$ for $i = 1,\dots,n$.

\item
(Grid of points)
Choose $\delta_n > 0$ and $\grid_n\subset\cube$ so that any $x\in\cube$ is the convex combination of points in $\grid_n$ whose distance to $x$ is not more than $\delta_n$.

\item
(Convexification)
Construct $\convn$ as in~\eqref{equ:convn:def}, for instance by using a convex hull procedure such as \qhull.
\end{enumerate}
\end{proc}

\begin{figure}
\capstart
\centering
\subfigure[Data and estimator $f_n$]{%
   \includegraphics[width=.48\textwidth]{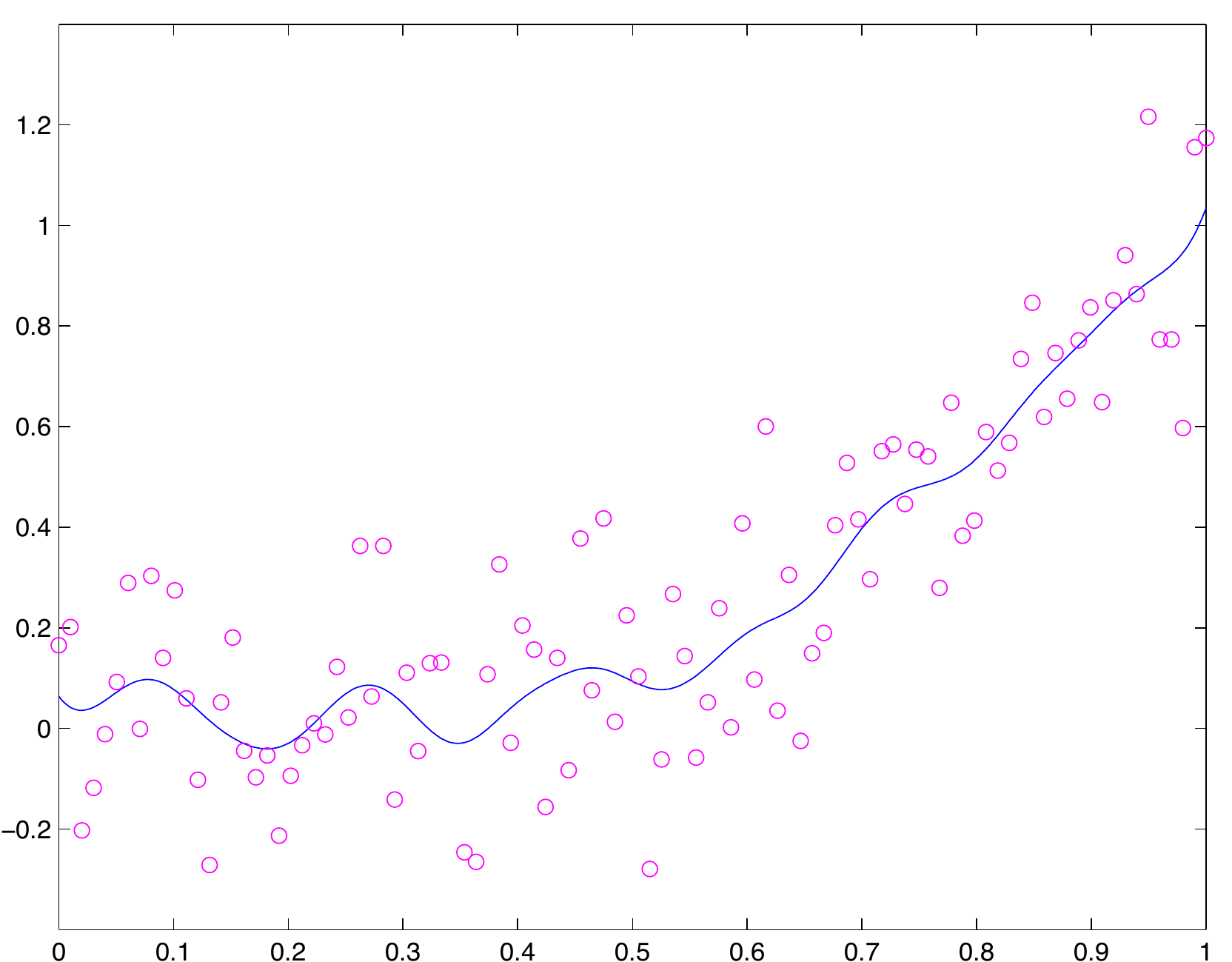}
   \label{fig:algorithm-step1}%
   }\hfill
\subfigure[Estimator $f_n$ and its values on $\grid_n$]{%
   \includegraphics[width=.48\textwidth]{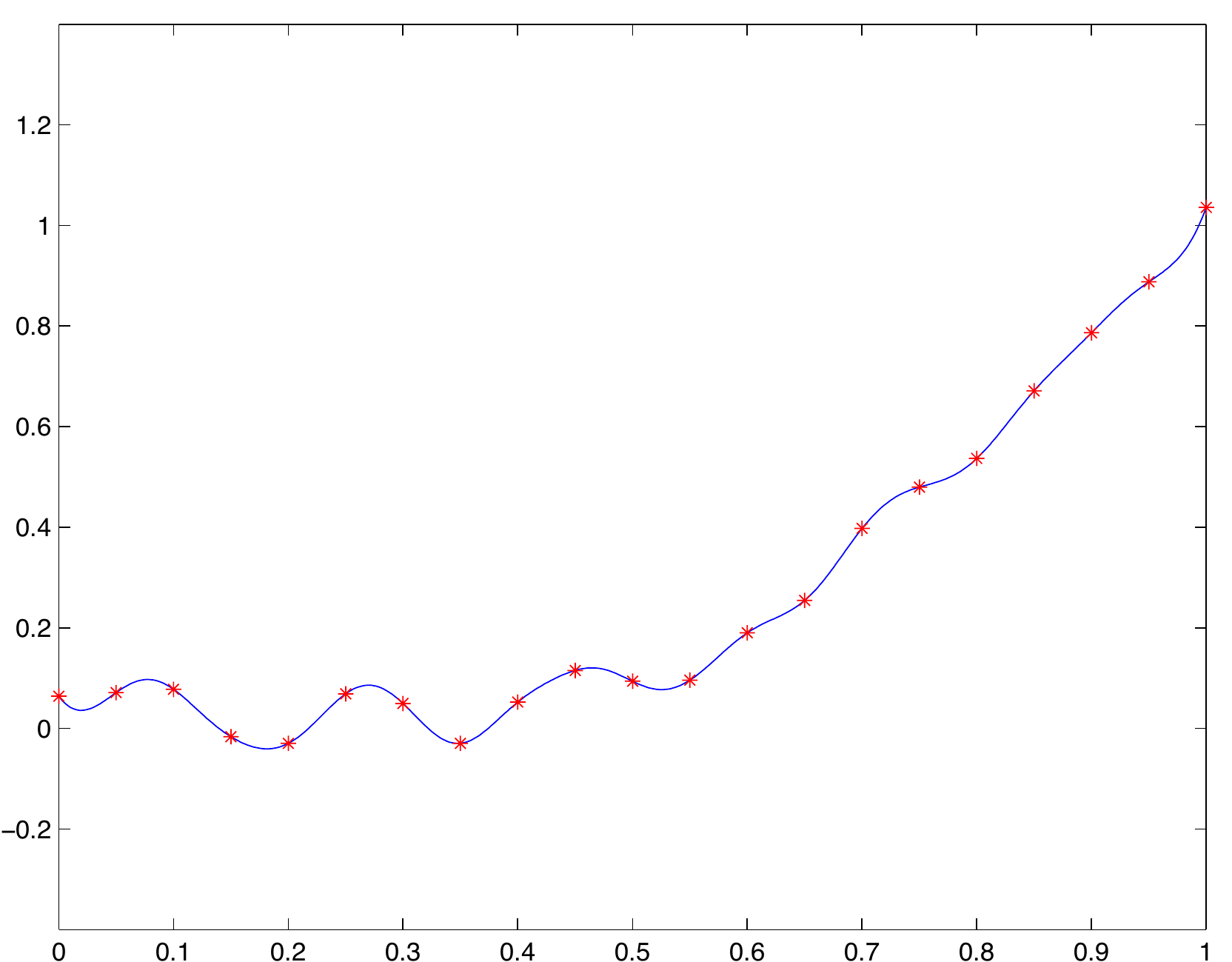}
   \label{fig:algorithm-step2}%
   }\\
\subfigure[Convex estimator $\convn$ from $f_n$ on $\grid_n$]{%
   \includegraphics[width=.48\textwidth]{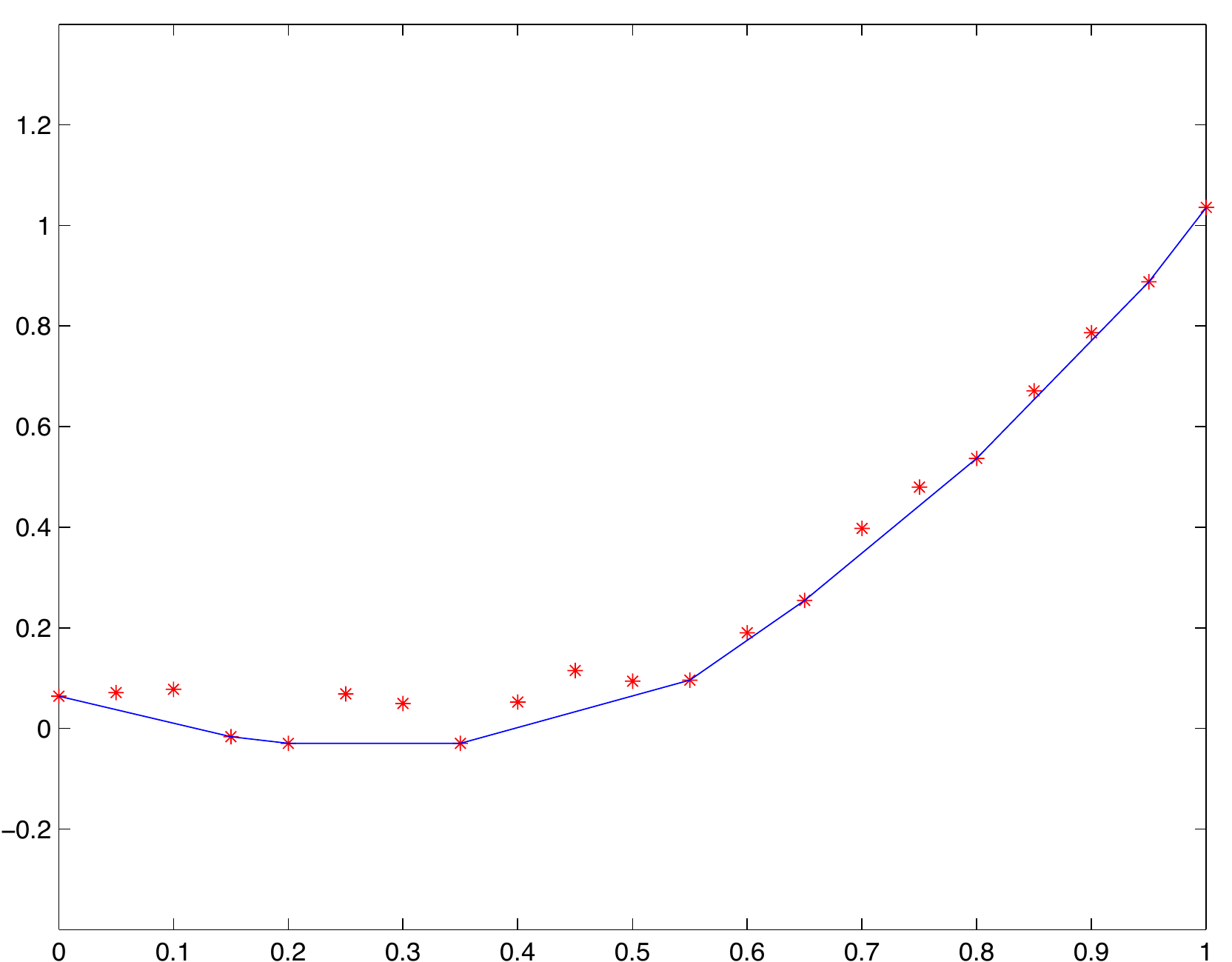}
   \label{fig:algorithm-step3}%
   }\hfill
\subfigure[Data and convex estimator $\convn$]{%
   \includegraphics[width=.48\textwidth]{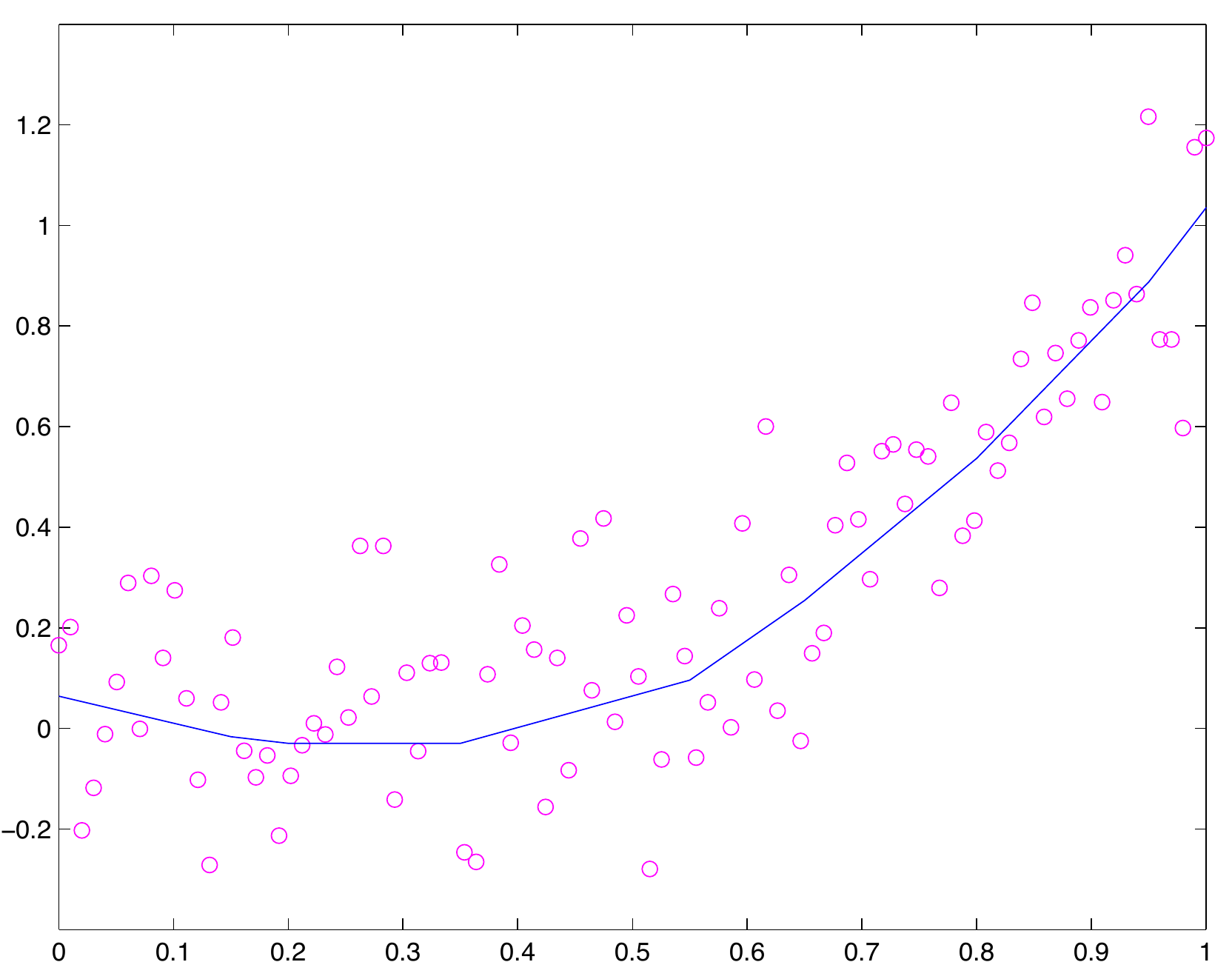}
   \label{fig:algorithm-step4}
   }
\caption{\small Steps in constructing a convex estimator}
\label{fig:algorithm}
\end{figure}

In \tref{Figure}{fig:algorithm} we represent the steps of the procedure with an example:
in~\subref{fig:algorithm-step1} we show the data and the resulting estimator $f_n$;
in~\subref{fig:algorithm-step2} we show the estimator and its values at the points of $\grid_n$;
in~\subref{fig:algorithm-step3} we show the convex estimator $\convn$ obtained from the values of $f_n$ at $\grid_n$;
and in~\subref{fig:algorithm-step4} we compare the original data and the convex estimator obtained.

We now show that if in the \tref{Procedure}{proc:convn}, $f_n$ is a good approximation of $f$, then $\convn$ is a good approximation of $f$ provided it satisfies:

\begin{cond}
\item\label{cond:f}
$f$ is a continuous convex function defined on $\cube$, with $\norm{f}_\Lip = L < \infty$, where
\[
   \norm{f}_\Lip
   = \sup \{\abs{f(x) - f(y)}/\abs{x - y} : x,y\in\cube,\,x\ne y \}.
\]
and $\abs{x - y}$ denotes the (Euclidean) distance between $x$ and $y$ in $\RR^d$.
(Recall that convex functions on $\cube$ are locally Lipschitz, but here we require that $f$ be uniformly Lipschitz in all of $\cube$.)
\end{cond}

\begin{thm}\label{thm:conv:1}
Suppose $f$ satisfies \condref{cond:f} and let $f_n$, $\delta_n$, $\grid_n$ and $\convn$ be as in \tref{Procedure}{proc:convn}, with
\begin{equation}\label{equ:thm:conv}
   \sup\,\{ \abs{f_n(x) - f(x)} : x\in\grid_n \} \le \varepsilon_n.
\end{equation}

Then,
\[
   - \varepsilon_n
   \le \convn(x) - f(x)
   \le \varepsilon_n + L \delta_n
   \quad\text{for all $x\in\cube$.}
\]
\end{thm}

\begin{proof}
Since $f$ is convex and $\varepsilon_n$ is a constant, the function $f - \varepsilon_n$ is convex.
Moreover, $f(x) - \varepsilon_n \le f_n(x)$ for all $x\in\grid_n$ implies that $f - \varepsilon_n\in\convl_n$, and by the definition of $\convn$ in~\eqref{equ:convn:def},
\[
   f(x) - \varepsilon_n \le \convn(x)
   \quad\text{for all $x\in\cube$,}
\]
proving one inequality.

For the other inequality, consider $x\in\cube$, and let $x_k\in\grid_n$ and $\lambda_k \ge 0$, $k = 1,\dots,d+1$, be such that
\[
   \sum_k \lambda_k x_k = x,
   \quad
   \sum_k \lambda_k = 1,
   \qand
   \abs{x - x_k}\le\delta_n\ \text{for $k = 1,\dots,d+1$.}
\]
Then,
\[
   \begin{aligned}
   \convn(x)
   &\le \sum_k \lambda_k\, \convn(x_k)
      && \text{since $\convn\in\convfs$,} \\
   &\le \sum_k \lambda_k\, f_n(x_k)
      && \text{by \tref{Lemma}{lem:convn:Ln},} \\
   &\le \sum_k \lambda_k\, (f(x_k) + \varepsilon_n)
      && \text{by~\eqref{equ:thm:conv},} \\
   &= \biggl(\sum_k \lambda_k\,f(x_k)\biggr) + \varepsilon_n
      && \text{since $\sum_k \lambda_k = 1$.}
   \end{aligned}
\]

Now, $\norm{f}_{\Lip} = L$ and $\abs{x_k - x} < \delta_n$, and therefore
\[ f(x_k) \le f(x) + L \delta_n. \]
Hence, since $\lambda_k \ge 0$ and using again that $\sum_k \lambda_k = 1$, we conclude
\[
   \convn(x)
   \le \biggl(\sum_k \lambda_k\,(f(x) + L \delta_n)\biggr) + \varepsilon_n
   = f(x) + L \delta_n + \varepsilon_n,
\]
and the result follows.\qedhere
\end{proof}

\begin{rem}
In the proof we have not used the finiteness of $\grid_n$, and only the values of $f_n$ on $\grid_n$ are used.
\end{rem}

Noticing that given $\delta_n > 0$ we may construct a finite set $\grid_n$ with the property that any $x\in\cube$ is a convex combination of points in $\grid_n$ whose distance to $x$ is no more than $\delta_n$, we have:

\begin{cor}\label{cor:conv:2}
If $f$ satisfies \condref{cond:f},
given an estimator $f_n$ of $f$ and $\delta_n > 0$,
we may find $\grid_n$ and define $\convn$
according to \tref{Procedure}{proc:convn}, so that
\[ \norm{\convn - f}_\infty \le \norm{f_n - f}_\infty + L \delta_n. \]
\end{cor}

\begin{rem}
In the extreme case where $f_n = f$ for all $n$, we have $\norm{f_n - f}_\infty = 0$, but $\norm{\convn - f}_\infty > 0$ in general (for instance, if $\grid_n$ is finite and $f$ is not piecewise linear).
\end{rem}

\tref{Corollary}{cor:conv:2} tells us that the convex estimator $\convn$ obtained through the \tref{Procedure}{proc:convn} inherits the approximation properties of the original estimator $f_n$, and the rate of convergence is preserved or even bettered provided $\delta_n$ is small enough.

To illustrate this behavior, let us consider the following well-known types of convergence of a sequence of nonnegative random variables ${(R_n)}_n$ to $0$, where ${(r_n)}_n$ is a bounded sequence of positive numbers (possibly converging to $0$), and we have denoted by $\PP$ the underlying probability measure:

\begin{tipo}
\item
\label{tipo:first}
\label{tipo:bounded} 
For every $\varepsilon > 0$ there exists $M > 0$ such that
$\sup_n\,\prob{R_n > M r_n} < \varepsilon$.

\item
\label{tipo:prob} 
$\lim_{n\to\infty} \prob{R_n> \varepsilon r_n} = 0$ for every $\varepsilon > 0$.

\item
\label{tipo:tran} %
\label{tipo:as} 
$R_n = O(r_n)$ or $R_n = o(r_n)$ \as

\item
\label{tipo:aco} 
For every $\varepsilon > 0$, $\sum_{n=1}^{\infty} \prob{R_n > \varepsilon r_n} < \infty$.
\savectr{enumi}
\label{tipo:last}
\end{tipo}

It is easy to see that:

\begin{thm}\label{thm:granfinale}
If any of \tiporef{tipo:first} through \tiporef{tipo:last} holds for $R_n = \norm{f_n - f}_\infty$, then it also holds for $R_n = \norm{\convn - f}_\infty$, provided $f$ satisfies \condref{cond:f} and $\convn$ is constructed as in \tref{Corollary}{cor:conv:2} with $\delta_n = o(r_n)$.
\end{thm}

For example, \citet{Tran} shows:

\begin{thm}
\label{thm:tran}
For $j=1, 2, \dotsc$, let $\{(\xx[j], \yy[j])\}_j$ be a strictly stationary sequence of random variables, where the $\xx[j]$ and the $\yy[j]$ are $\RR^d$-valued and $\RR$-valued, respectively.
Suppose $f(x) = \Esp{\yy\mid\xx=x}$ is estimated by
\[
   f_n(x)=\frac{1}{\card{I_n(x)}} \sum_{i\in I_n(x)} \yy[i]
   \quad\text{for $x \in\cube$},
\]
where
$I_n(x) = \{ i : 1\le i \le n, \abs{\xx[i] - x} \le \hn\}$,
and
$\hn \approx \left( \log(n) / n \right)^{1/(d+2)}$.

Then, under appropriate assumptions (including adequate regularity conditions),
\[ \norm{f_n - f}_{L^\infty(\cube)} = O(\hn) \quad \text{\as} \]
\end{thm}

Tran's result gives a \tiporef{tipo:tran} type of convergence, and therefore (by \tref{Theorem}{thm:granfinale}) we have that under the same assumptions,
\[ \norm{\convn - f}_{L^\infty(\cube)} = O(\hn) \quad \text{\as}, \]
provided we take $\delta_n = o(\hn)$ in \tref{Corollary}{cor:conv:2}.

More elaborate types of convergence include exact asymptotic behavior.
A very simple model might be, assuming $\xx[n]$ uniformly distributed on $\cube$:

\begin{tipo}
\backctr{enumi}
\item\label{tipo:dist}
There exist a sequence ${(d_n)}_n$ converging to $0$, and a random variable $R$ such that
\[ \prob{ r_n^{-1} (R_n - d_n)\le t } \to \prob{R \le t}, \]
for every $t\in\RR$ at which $\prob{R\le t}$ is continuous.
\end{tipo}

It is not possible in general to carry over this convergence from $R_n = \norm{f_n - f}_\infty$ directly to $R_n = \norm{\convn - f}_\infty$, as in general $\norm{\convn - f}_\infty$ could be much smaller than $\norm{f_n - f}_\infty$, and we cannot control $\norm{f_n - f}_\infty$ solely in terms of $\norm{\convn - f}_\infty$ and $\norm{f}_\Lip$.
Needless to say, by enlargening $r_n$ we may transform a \tiporef{tipo:dist} type into, say, a \tiporef{tipo:prob} type of convergence.

Besides the interest in itself, the convergence of type \tiporef{tipo:dist} allows us to find uniform confidence bands for the regression curve, which is a practical concern.
More precisely, if \tiporef{tipo:dist} is verified, for any $\alpha$, $0 < \alpha < 1$, we may find optimal (or near optimal) $s$ so that
\begin{equation}\label{equ:band}
   \prob{R_n \le s} \ge 1 - \alpha.
\end{equation}

If this inequality holds for $R_n = \norm{f_n - f}_\infty$ and assuming $\convn$ is constructed as in \tref{Corollary}{cor:conv:2} with $\delta_n = o(1)$ for all $n$, then~\eqref{equ:band} is valid for $R_n = \norm{\convn - f}_\infty$, albeit not with optimal $s$.

In other words, \tref{Corollary}{cor:conv:2} allows us to convert a uniform confidence band for $f_n$ of the form~\eqref{equ:band} into a (slightly different) uniform confidence band for $\convn$.

For instance, \citet[Theorem~2.1]{Jh82} shows:

\begin{thm}
\label{thm:johnston}
Let $(\xx[1],\yy[1]),\dots,(\xx[n],\yy[n])$ be a random sample from a bivariate population, with $\xx$ uniformly distributed in $\cube = [0,1]$, and consider the following estimator of $f(x) = \esp{\yy \mid \xx = x}$,
\begin{equation}\label{equ:johnston}
f_n(x) = \frac{1}{n \hn} \sum_{i=1}^n \yy[i]\, K((x-\xx[i])/\hn),
\end{equation}
where
$\hn \approx n^{-\delta}$
for some $\delta$, $1/5 < \delta < 1/3$,
and $K$ is a piecewise smooth density function with support in $[-A,A]$, $A > 1$.

Then,
under appropriate regularity assumptions we have
\[
   \Prob{ (2\delta \log n)^{1/2}
   \left[
      \sup_{0\le x\le 1} r_n^{-1}(x)\, \big(f_n(x)-f(x)\big) - d_n \right]
      < t}
   \to e^{-2\exp{(-t)}},
\]
where
\begin{equation}\label{equ:johnston:2}
   r_n^2(x) =
   \frac{\int K^2(u)\,\dd{u}\times\esp{\yy^2\mid\xx = x}}{n\hn}
\end{equation}
and
$d_n = O\left((2\delta\log n)^{1/2}\right)$.
\end{thm}

Confidence bands follow immediately \citep[Corollary 3.1]{Jh82}:

\begin{cor}\label{cor:johnston}
Assuming \tref{Theorem}{thm:johnston} holds, an approximate $(1-\alpha)\times 100 \%$ confidence band is
\[
   f_n(x) \pm r_n\,\bigl(d_n + c(\alpha) (2 \delta \log n)^{-1/2}\bigr),
\]
where $c(\alpha) = \log 2 - \log \abs{\log (1-\alpha)}$ (for practical applications, one would estimate $\esp{\yy^2\mid \xx = x}$ in \eqref{equ:johnston:2}).
\end{cor}

\tref{Theorem}{thm:johnston} and its corollary are still valid if instead of \eqref{equ:johnston}, $f_n$ is a $M-$smoother estimator defined as a solution of
\[ 0 = \frac{1}{n \hn} \sum_{i=1}^n \psi(\yy[i]-f_n)\, K((x-\xx[i])/\hn),\]
with $\psi$ a bounded monotone, antisymmetric real function \citep{Ha89}.

As a final remark, let us point out that we have only used that $f_n$ approximates the Lipschitz convex function $f$, independently of whether $f_n$ has been obtained through a smoothing procedure or any other approximation method.

\section{Numerical results}
\label{sec:numerical}

In this section we report on some practical aspects of our algorithm and present some simulations and examples showing its performance.

\subsection{Implementation}

We implemented our algorithm using \matlab.
The smoothing step was done with local polynomials of degree 1 with Gauss's kernel, and for the convexification we used \matlab's functions \rutina{convhull} (dimension 1) and \rutina{convhulln} (higher dimensions), which are based upon the \qhull algorithm described in~\citet{qhull}.

The bandwidth was chosen using cross-validation for the local-polynomial fitting at the data points.
In the examples shown below, once the optimal bandwidth was chosen, the local-polynomial fitting function was computed at the same data points which were set a priori as design.
Whenever the data points were not a priori designed, the local-polynomial fit was evaluated on a uniform grid having approximately the same number of points.

\subsection{One dimensional simulations}

In this section we briefly illustrate the finite sample properties of the convex estimate of the regression function by means of a simulation study.
For this purpose we considered the same three examples presented in~\citet{BD07}, namely,
\begin{align*}
   f_1(x) &= e^{3(x-1)}, \\
   f_2(x) &= \frac{16}{9} \left(x-\frac14\right)^2, \\
   f_3(x) &= \begin{cases}
      -4x+1 & \text{if } 0 \le x \le 1/4, \\
      0 & \text{if } 1/4 < x < 3/4, \\
      4x-3 & \text{if } 3/4 \le x,
      \end{cases}
\end{align*}
and $Q = [0,1]$.
Notice that even though the third function is just Lipschitz, all these functions satisfy the assumption \condref{cond:f}.

As in~\citet{BD07}, we ran some simulations with $n = 100$ uniformly distributed design points for the explanatory variables and added a normal noise with standard deviation $\sigma=0.1$ to the response variable.

\begin{figure}
\capstart
\begin{center}
\includegraphics[width=.32\textwidth]{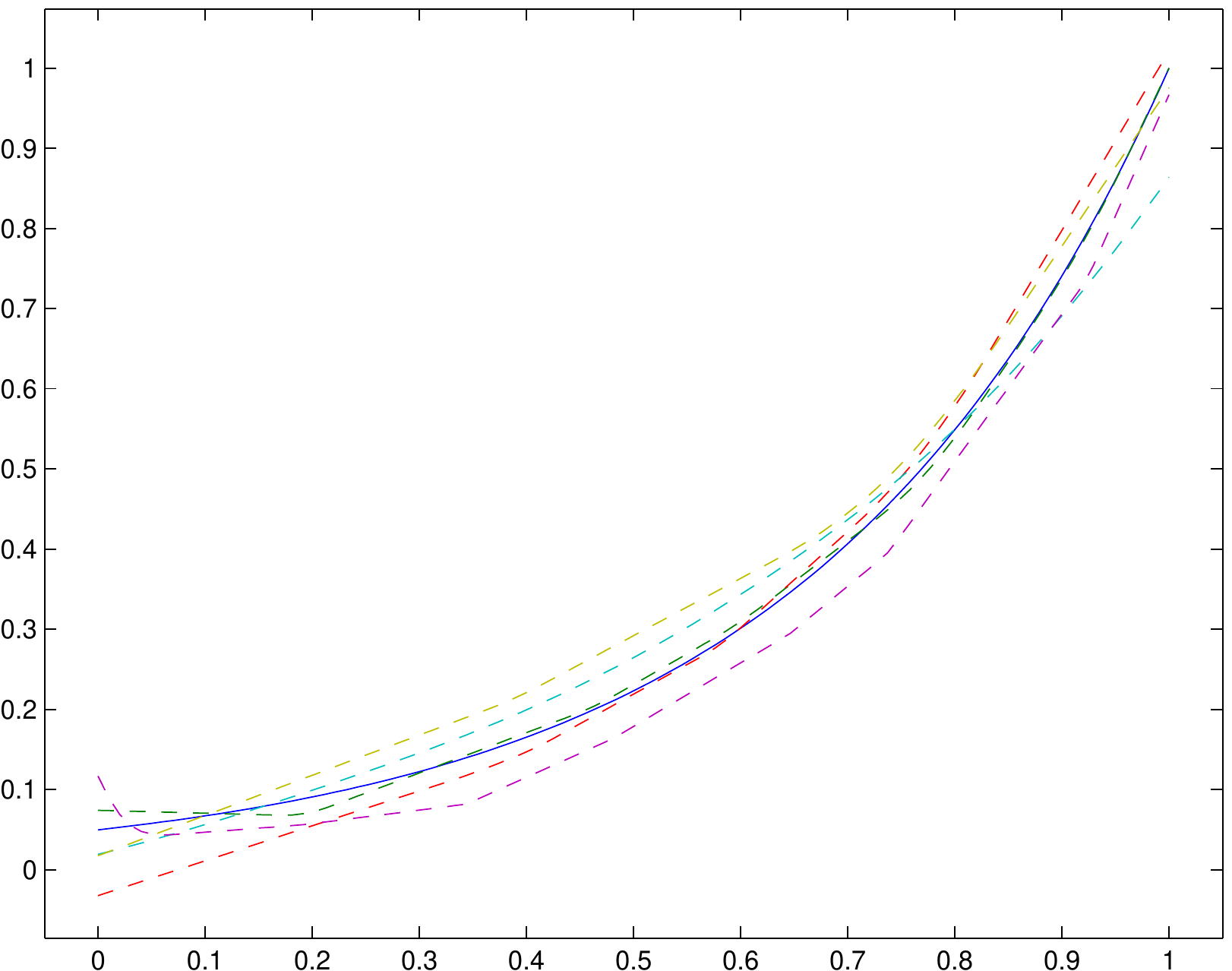}
\includegraphics[width=.32\textwidth]{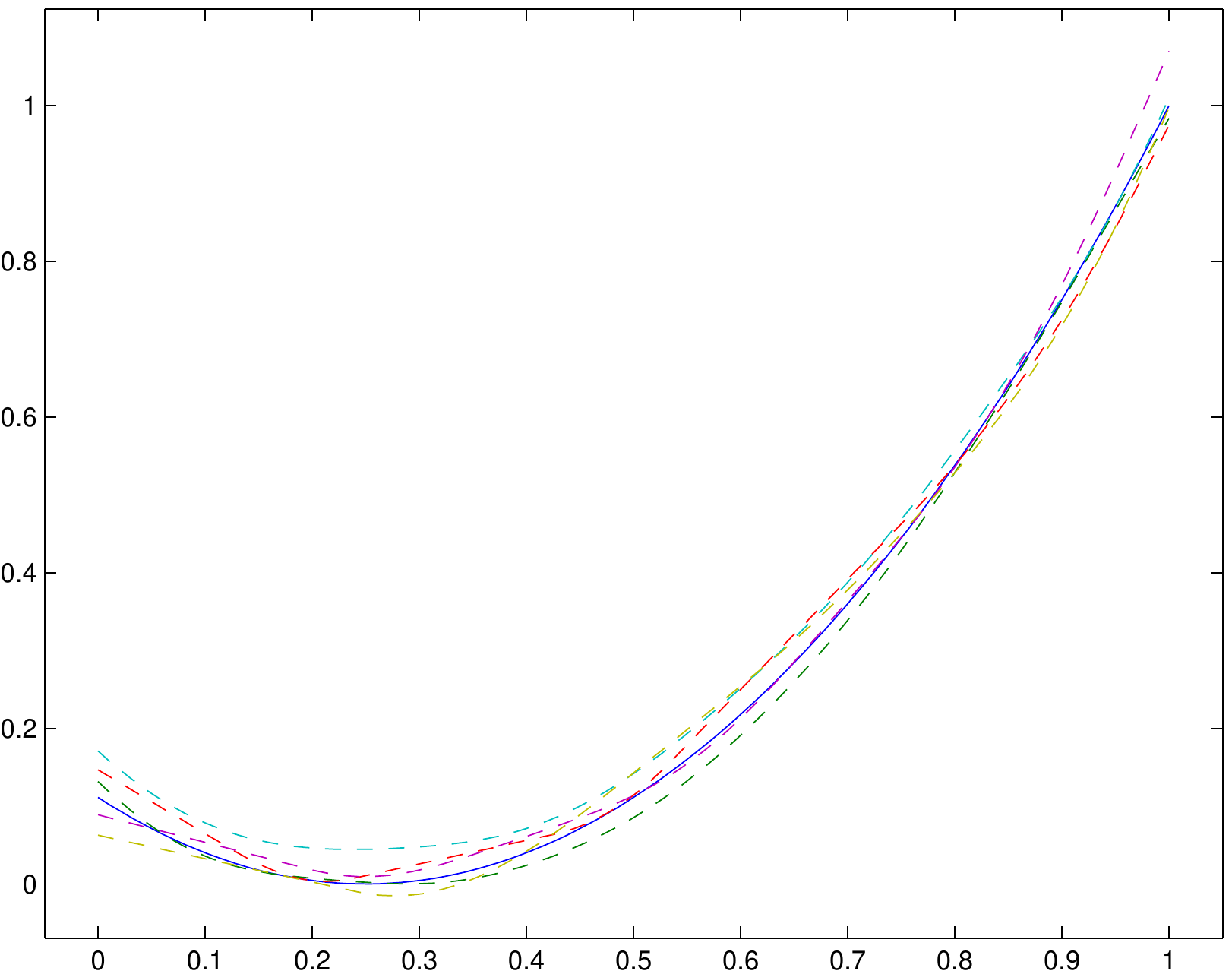}
\includegraphics[width=.32\textwidth]{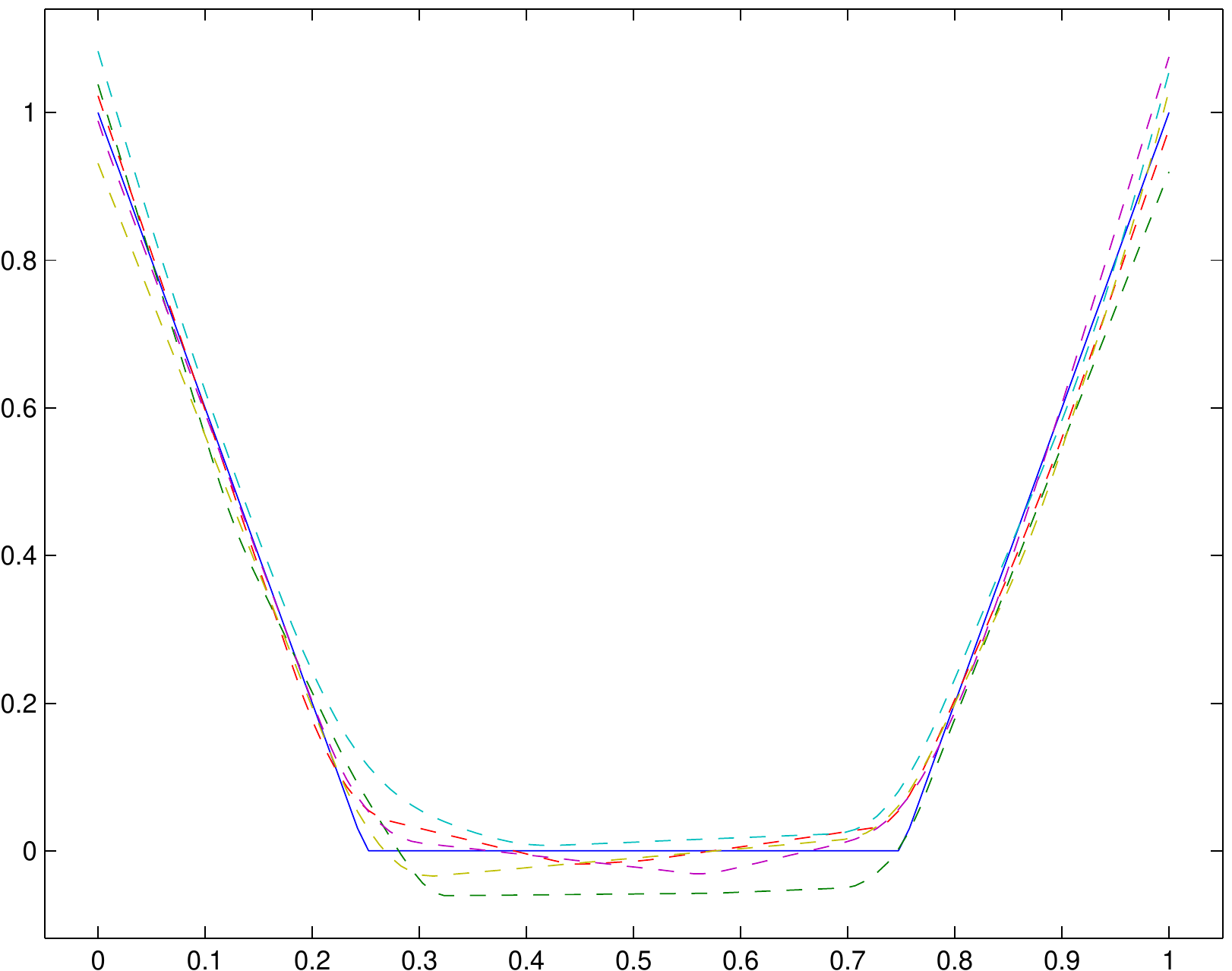}
\end{center}
\caption{\small
Regression functions $f_1$ (left), $f_2$ (middle), $f_3$ (right), and their estimates.
Result of 5 simulations for each regression function, with sample size $n = 100$ and normal errors with $\sigma = 0.1$.
The estimates are very reasonable, even for $f_3$, which is just Lipschitz, and not $C^1$.}
\label{F:regression1d}
\end{figure}

In \tref{Figure}{F:regression1d} we display for each regression function five typical estimates obtained from different simulation runs observing a typical performance.
The estimates for the two smooth functions $f_1$ and $f_2$ are comparable to the regressions obtained in~\citet{BD07}, but our estimates of the nonsmooth regression function $f_3$ exhibit a much closer fit.
This is an advantage of our method, which does not approximate the derivative of the regression function, and thus it demands less smoothness and approximates better non differentiable functions.

\begin{figure}
\capstart
\begin{center}
Var \hfil\hfil Bias$^2$ \hfil\hfil MSE

\includegraphics[width=.32\textwidth]{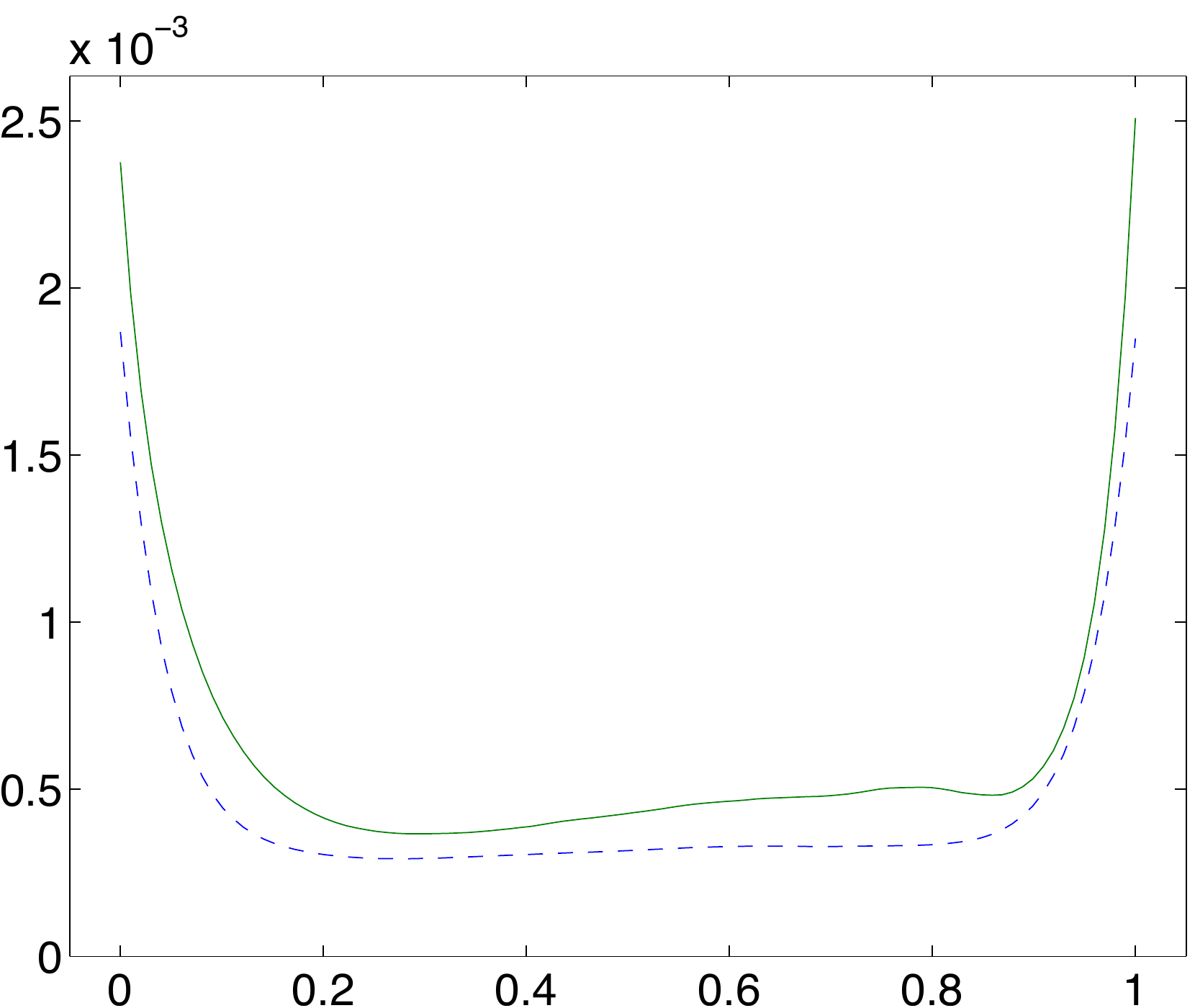}
\includegraphics[width=.32\textwidth]{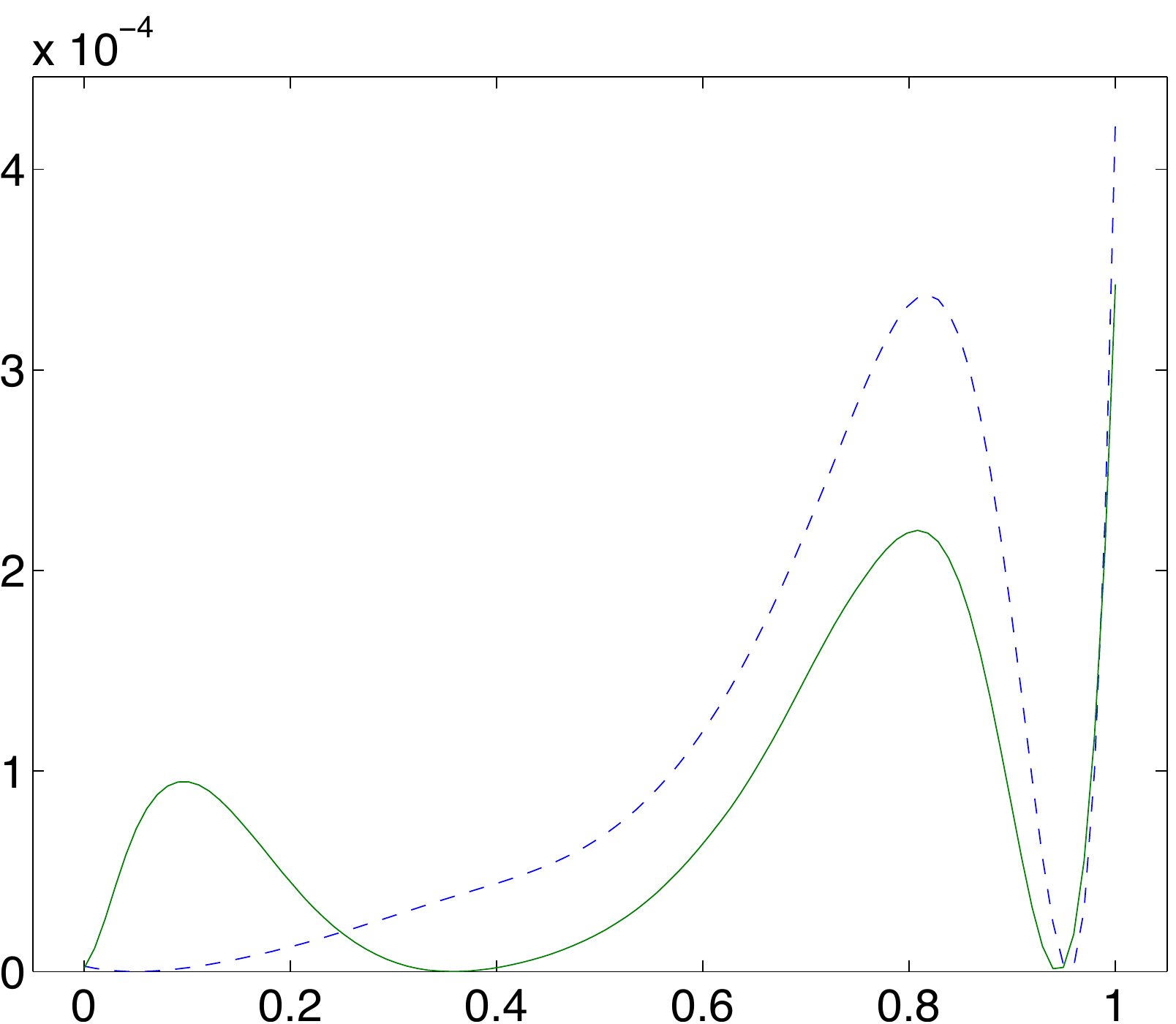}
\includegraphics[width=.32\textwidth]{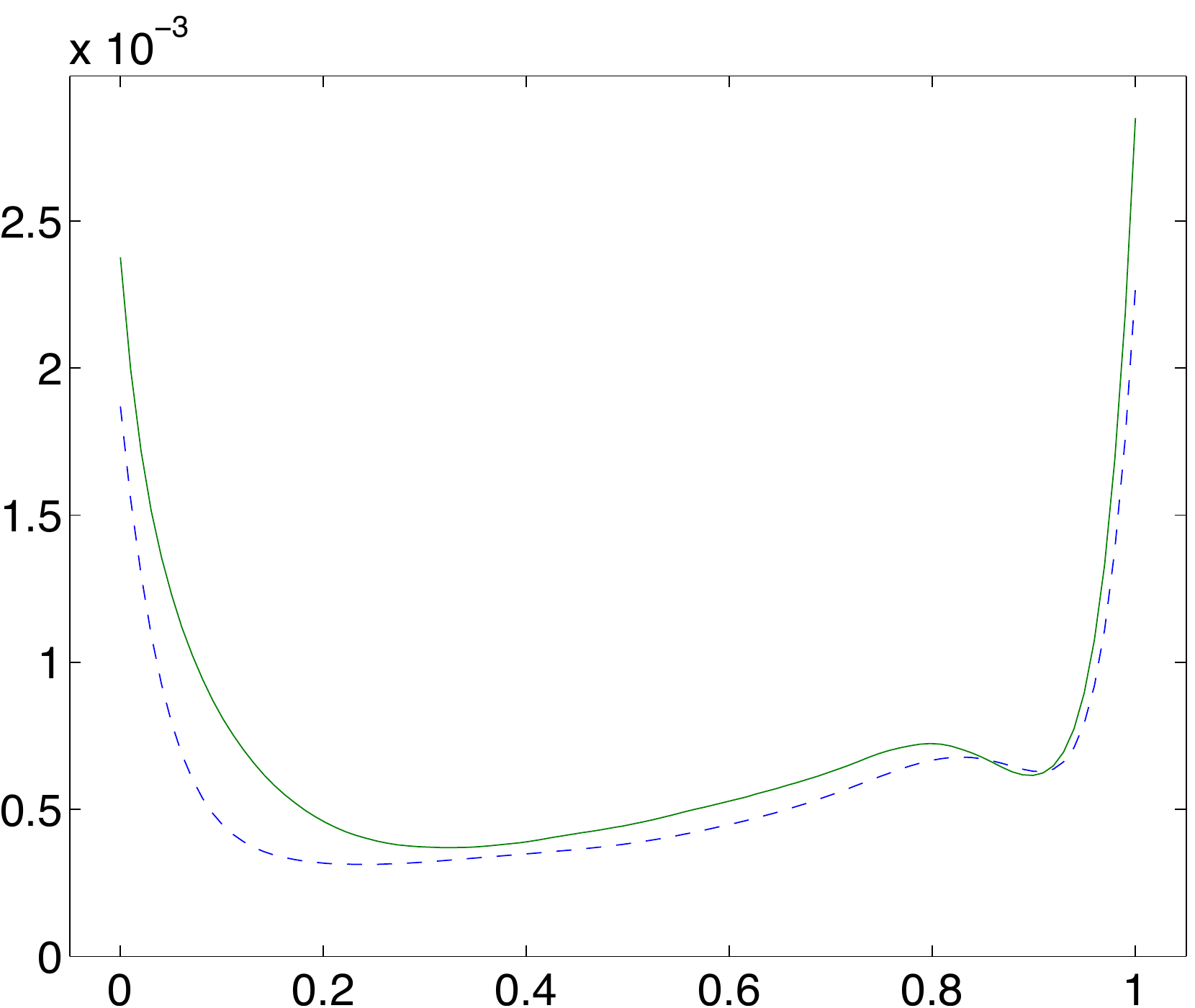}

\includegraphics[width=.32\textwidth]{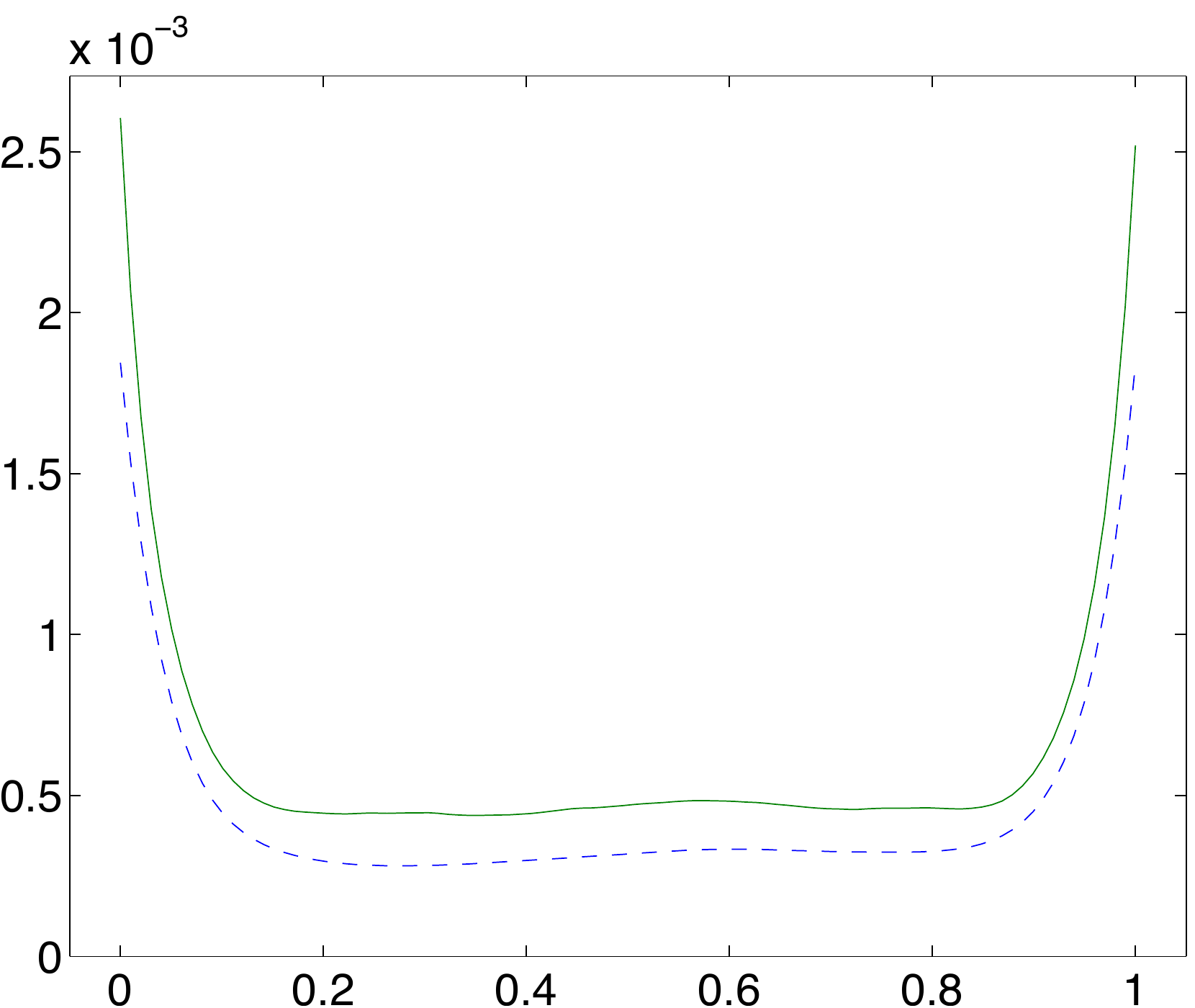}
\includegraphics[width=.32\textwidth]{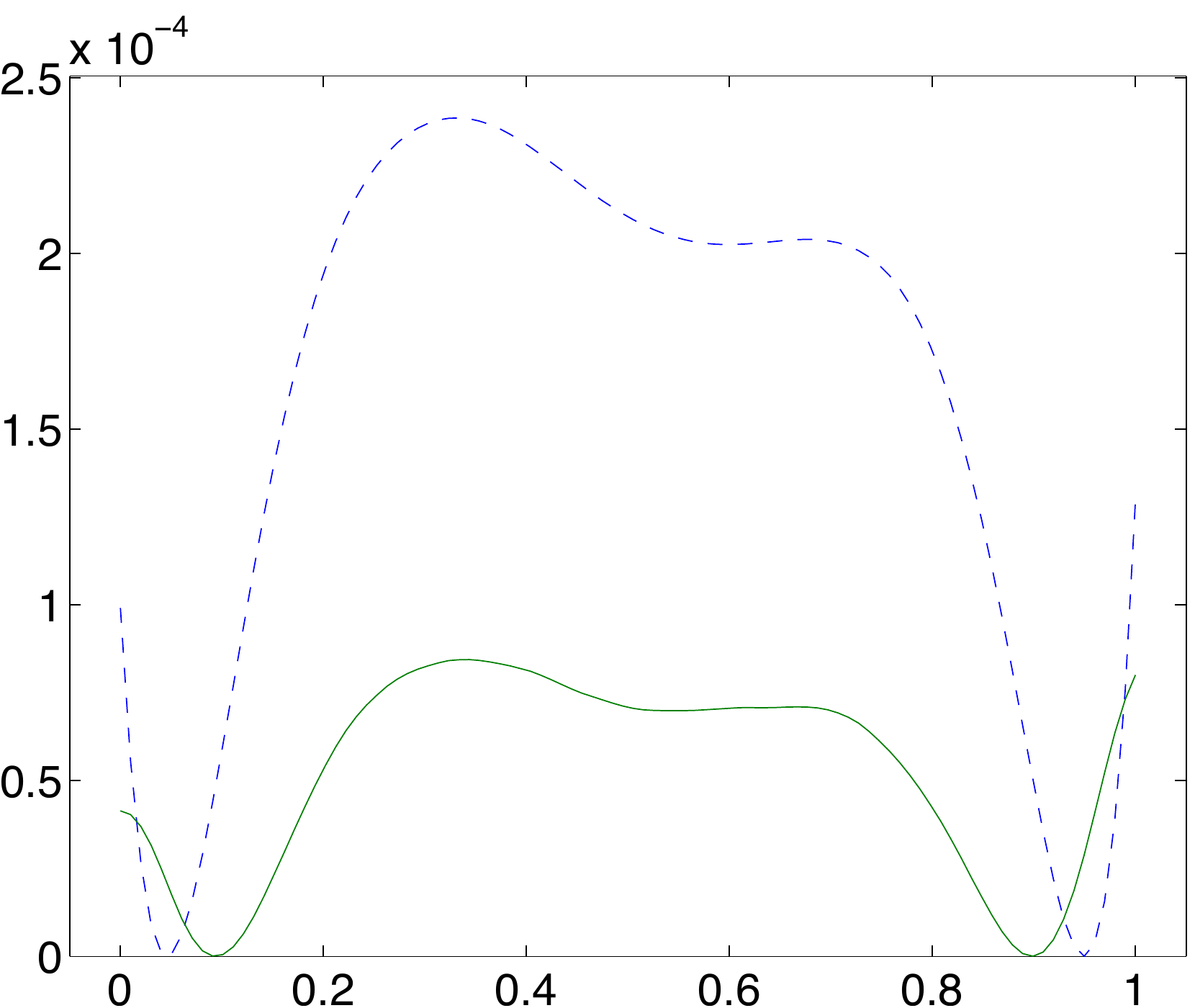}
\includegraphics[width=.32\textwidth]{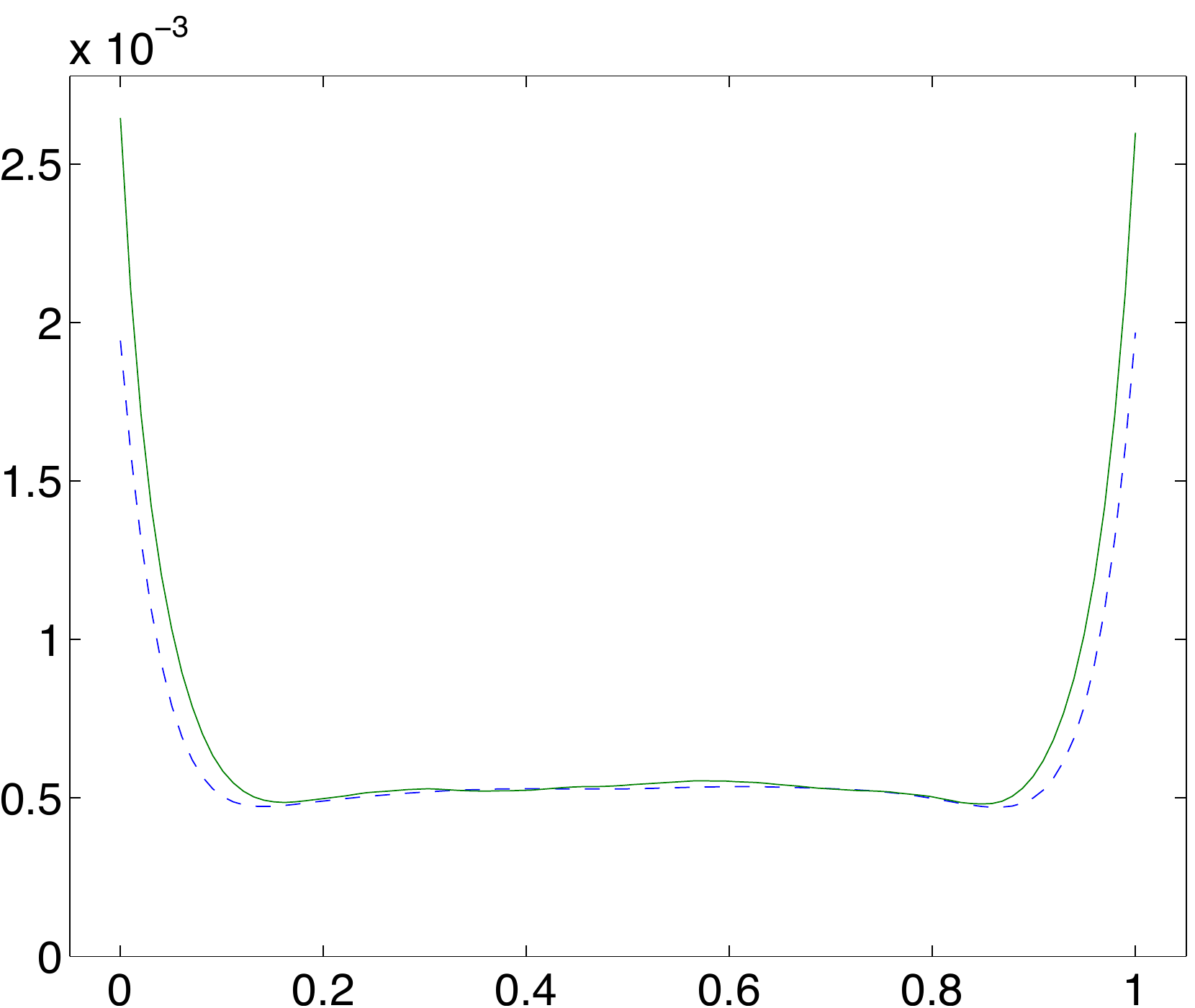}

\includegraphics[width=.32\textwidth]{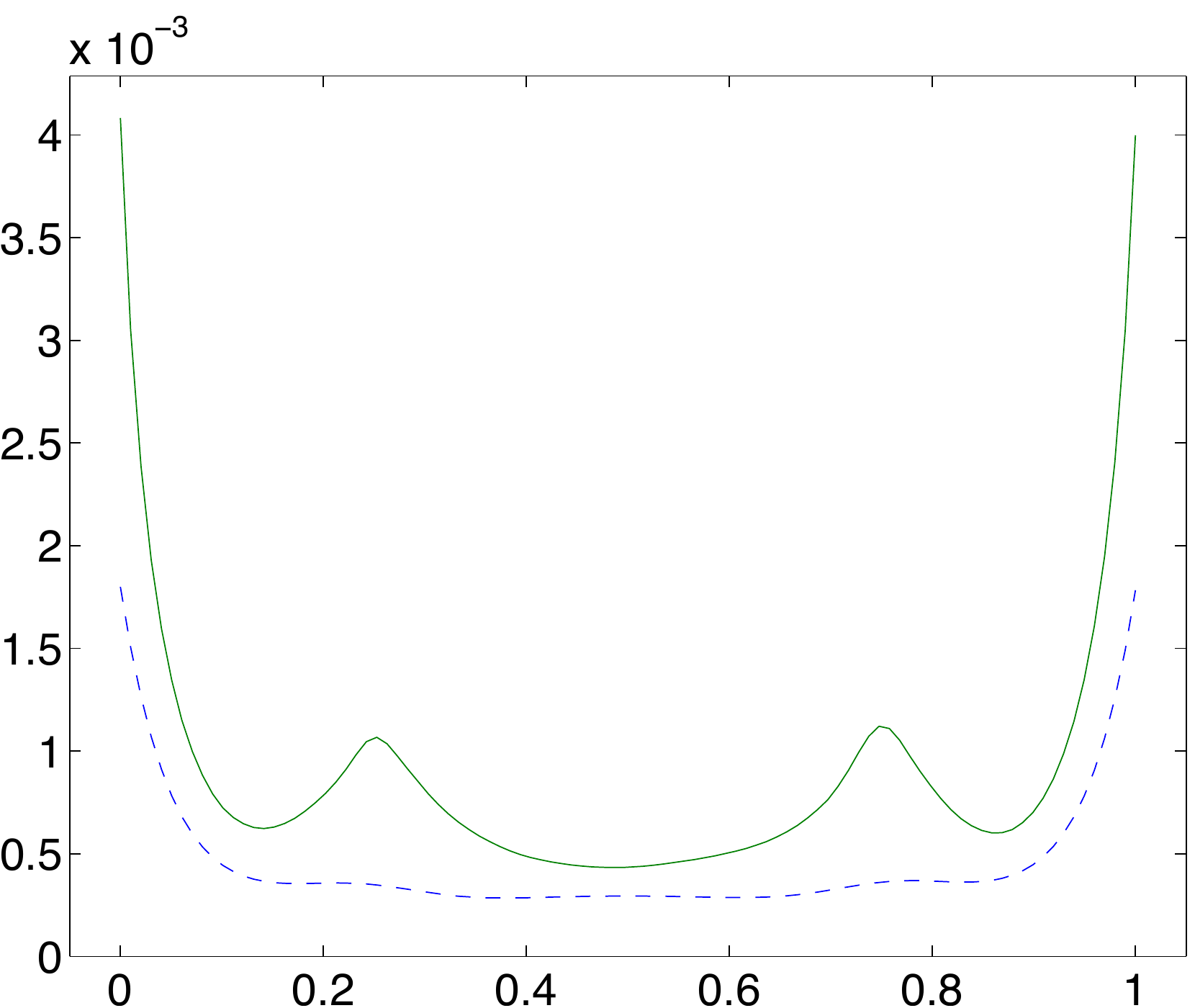}
\includegraphics[width=.32\textwidth]{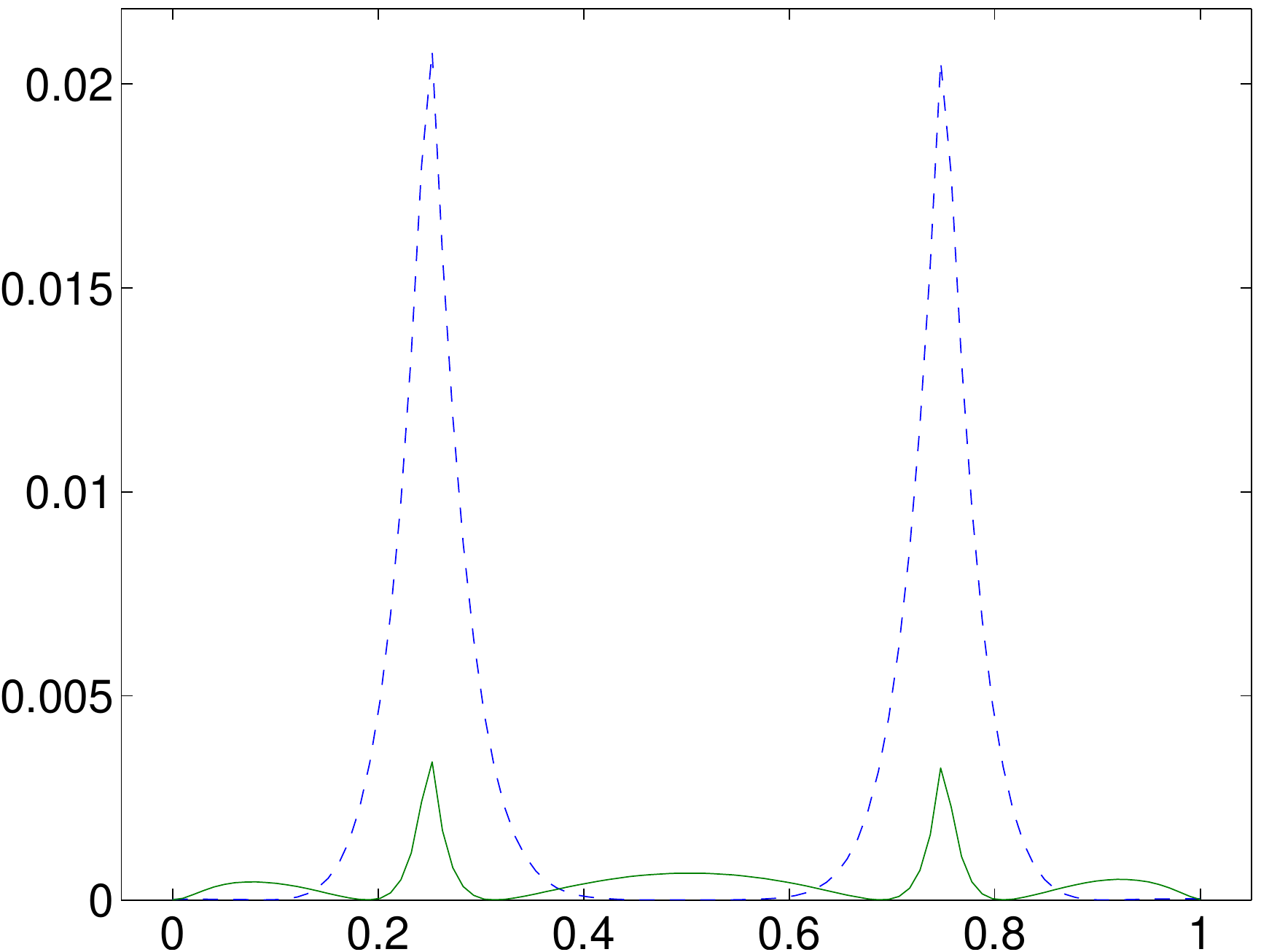}
\includegraphics[width=.32\textwidth]{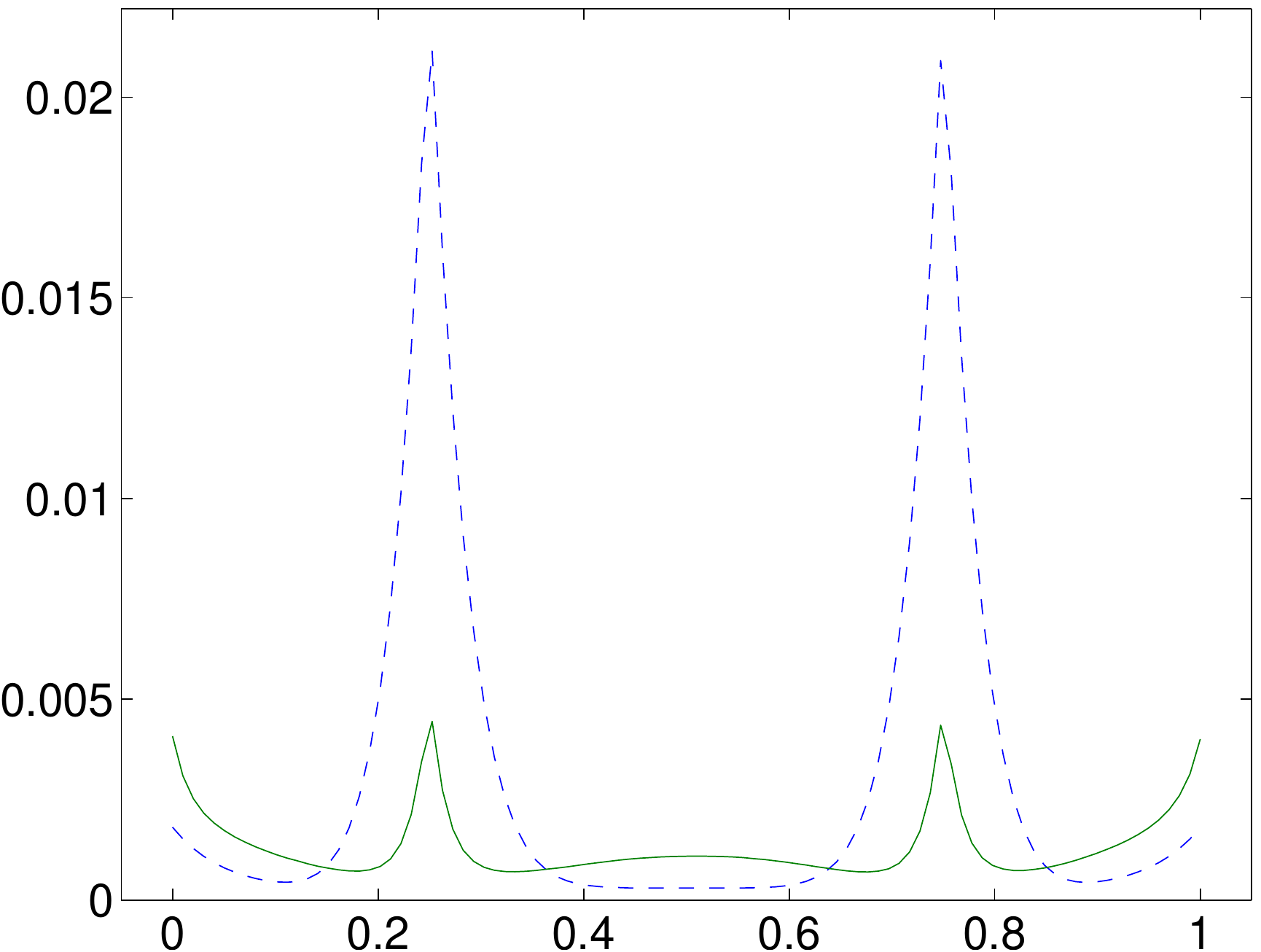}
\end{center}
\caption{\small
Variance (left), squared Bias (middle) and Mean Square Error (right) of our convex estimate (solid line) and of the local linear estimate (dashed line).
These indicators were obtained with 2000 simulation runs, for $f_1$ (top), $f_2$ (middle) and $f_3$ (bottom), using 100 uniformly distributed design points for the explanatory variable, and normal error with $\sigma = 0.1$ for the observed variable.
Only small differences between the local linear estimate and our convex estimate are observed.
In some cases, our convex estimate is even better
.}
\label{F:varbiasmse1d}
\end{figure}

In the second part of this simulation study we investigated the mean square error, bias and variance of our convex estimate.
For this we considered again the three regression functions in $f_1$, $f_2$, $f_3$ and computed---with 2000 simulation runs---the curves for the mean square error, squared bias and variance.
The results shown in \tref{Figure}{F:varbiasmse1d} look very much alike those in~\citet{BD07}, except for the ones related to $f_3$, where our estimator seems to be better.
In this figure the mean square error, bias and variance of the estimator by local linear polynomials are represented by the dashed lines, while those quantities related to our convex estimator are represented by the solid lines.

Finally, in \tref{Figure}{F:confidence} we show approximate 95\% confidence bands for one estimate to each of the previous regression functions.
We ran a simulation with 100 uniformly distributed design points for the explanatory variables and added a normal noise with $\sigma=0.1$ to the response variable.
In order to use the existing results on the width of the confidence bands from \citet[Corollary 3.1]{Jh82} (see also \tref{Theorem}{thm:johnston} and its corollary), the smoothing step was done with the formula
\[ f_n(x) = \frac{1}{n \hn} \sum_{j=1}^n K((x-\xx[i])/\hn), \]
where $K(x) = \frac34 (1-x^2)_+$ is Epanechnikov's Kernel.
This regression formula has bad approximation properties at the endpoints of the interval, which explains the mild misfit observed there.
The width of the band was $0.1392$, $0.1382$, $0.1628$, for the estimate corresponding to $f_1$, $f_2$, and $f_3$, respectively.

\begin{figure}
\capstart
\begin{center}
\includegraphics[width=.32\textwidth]{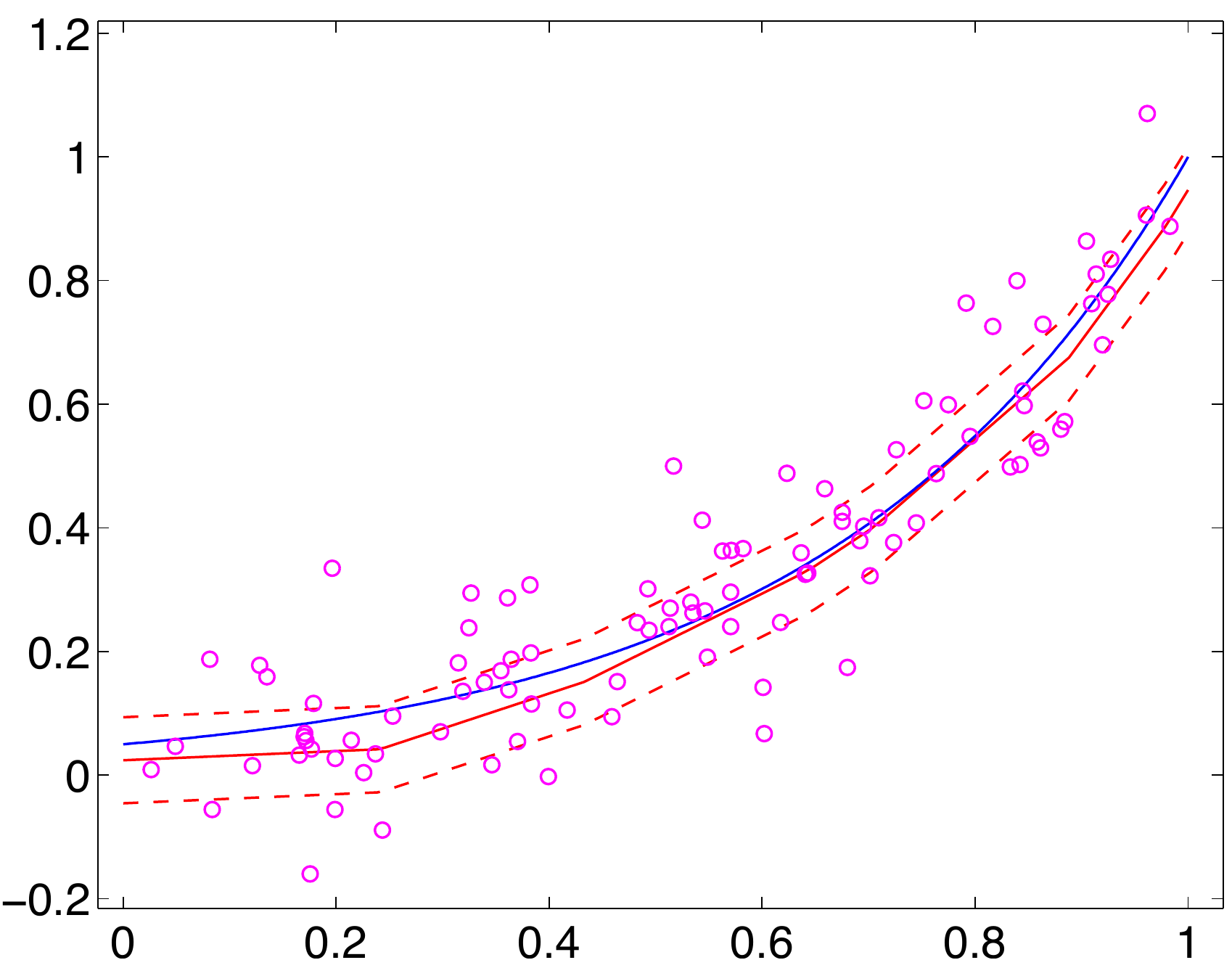}
\includegraphics[width=.32\textwidth]{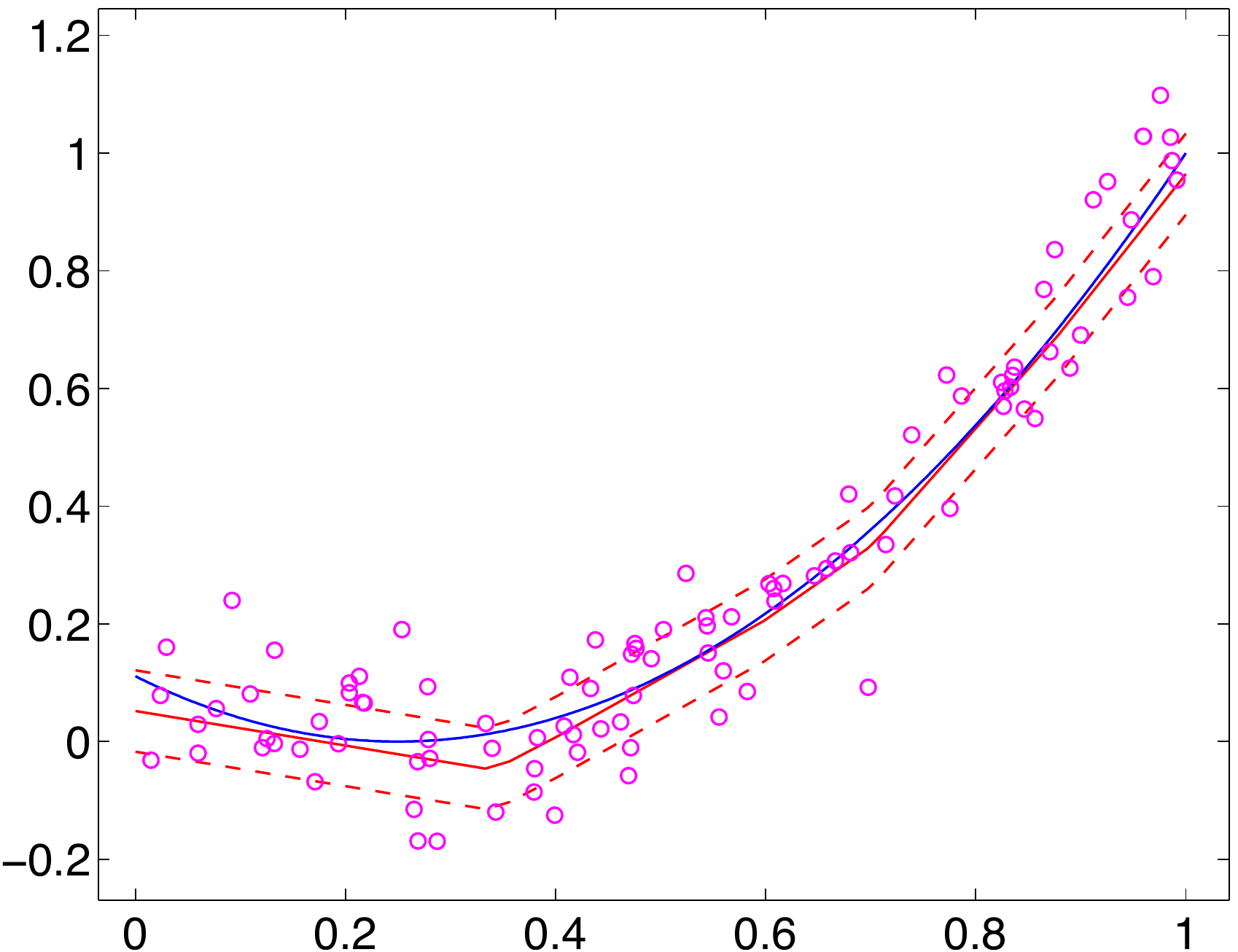}
\includegraphics[width=.32\textwidth]{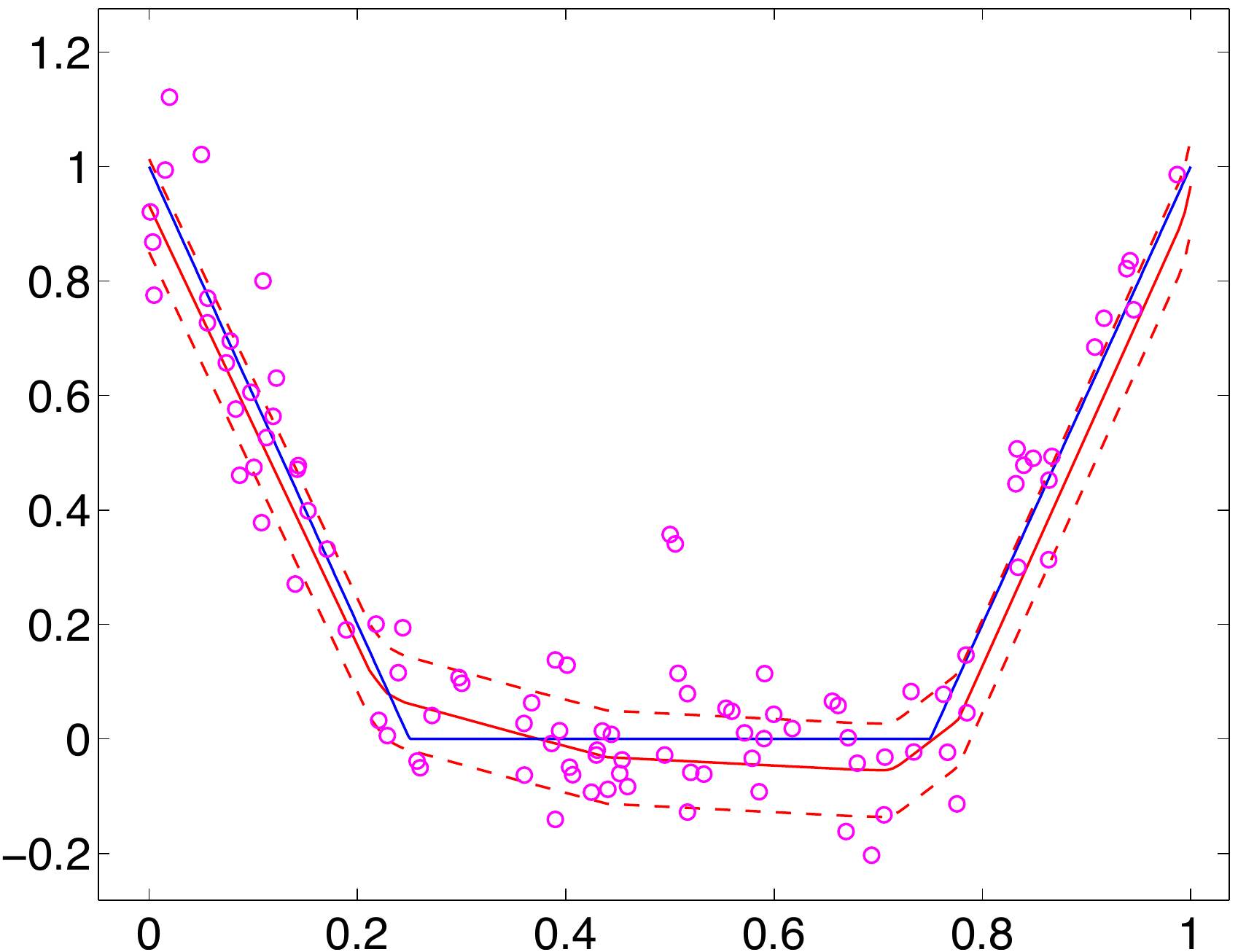}
\end{center}
\caption{\small
Approximate 95\% confidence bands.
In each plot we show the exact regression function, the estimate, the 95\% confidence bands, and the data points for $N=100$ design points, and normal error with $\sigma=0.1$.
The bands have width $0.1392$, $0.1382$, and $0.1628$ for $f_1$ (left), $f_2$ (middle), $f_3$ (right), which were computed with the formula provided in \tref{Corollary}{cor:johnston}.}
\label{F:confidence}
\end{figure}

\subsection{Rabbits' data}

We studied an example considered in~\citet{DM61}, who analyzed the relationship between age and eye lens weight for rabbits in Australia.
This relationship is expected to be guided by a concave function.
In this study, the dry weight of the eye lens was measured (in milligrams) for 71 free-living wild rabbits of known age (measured in days).
A detailed description of the experiment and the data can be found in \url{http://www.statsci.org/data/oz/rabbit.html}.
The data was analyzed by~\citet{Ratkowsky83} using a parametric nonlinear growth model, and by~\citet{BD07} with their non-parametric convex regression method.
We used our method to obtain the concave regression, with the smoothing step performed with local polynomials of degree 1 and 2, and report the findings in \tref{Figure}{F:rabbits}.
In both cases, the bandwidth for the local polynomial smoothing was set using cross-validation, and the result of the smoothing step was evaluated at a uniform grid of 100 points.
The convexification step yielded the estimated regression curves that can be observed in \tref{Figure}{F:rabbits} with an excellent fit to the data.

\begin{figure}
\capstart
\begin{center}
\includegraphics[width=.49\textwidth]{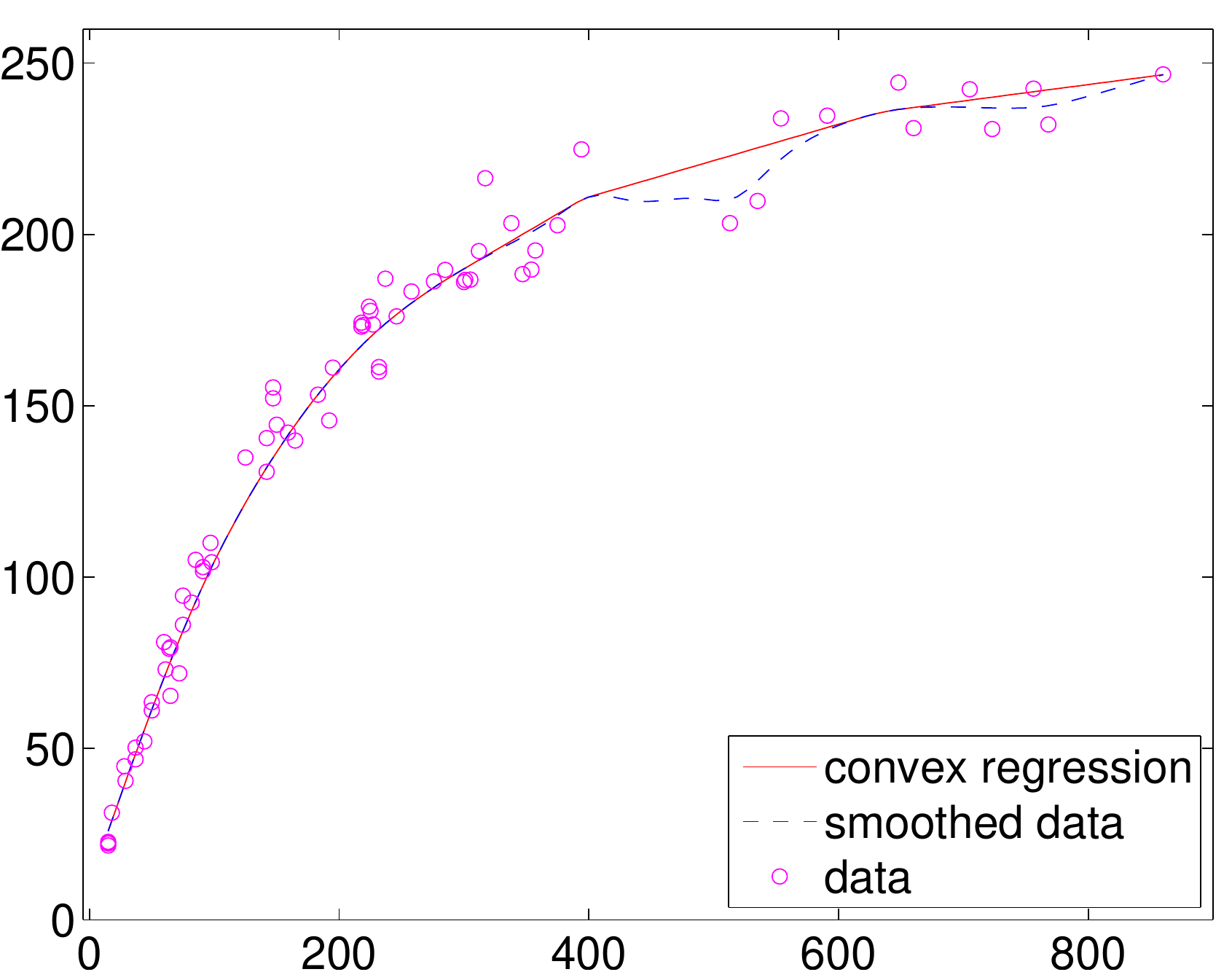}
\includegraphics[width=.49\textwidth]{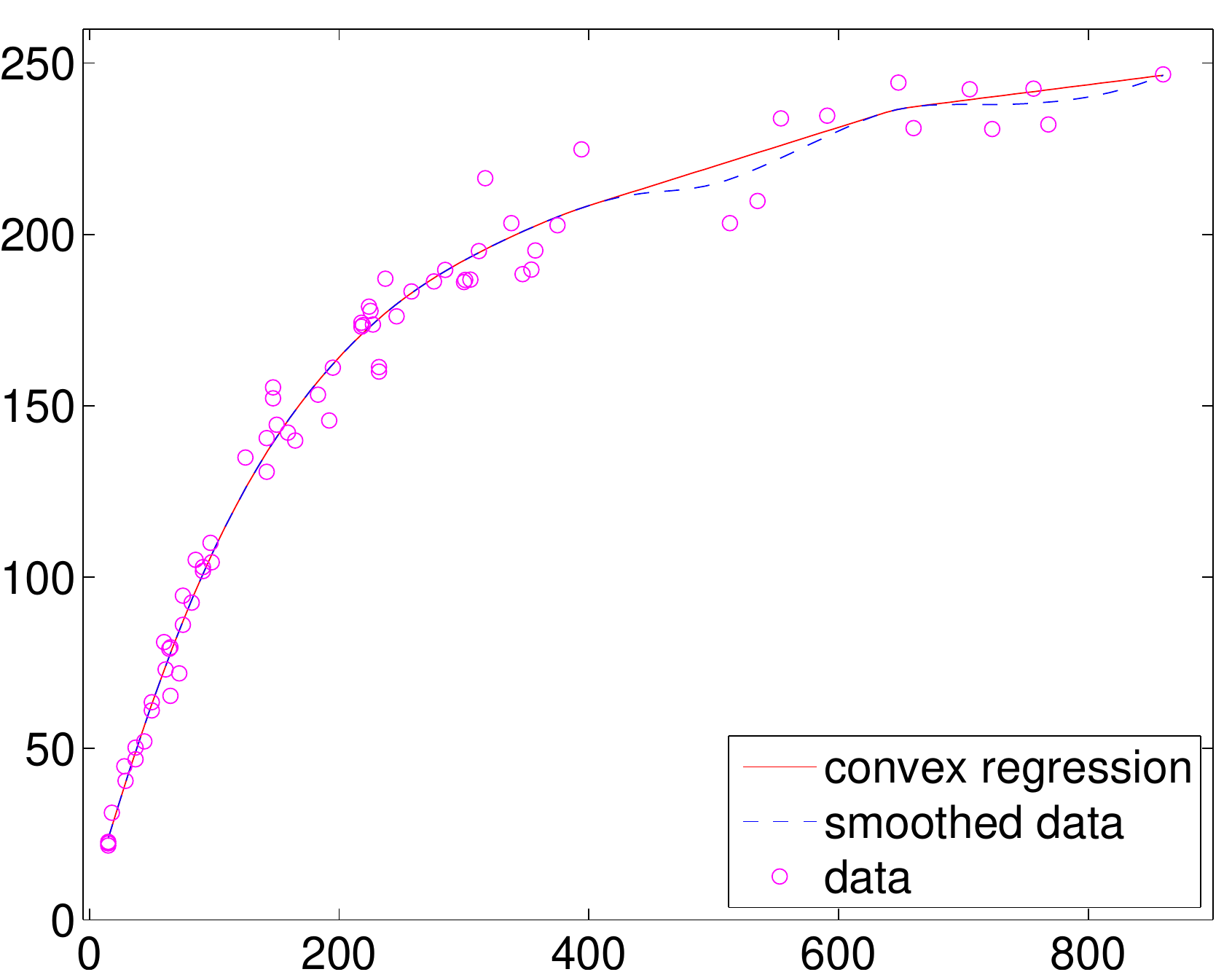}
\end{center}
\caption{\small
Convex regression of the rabbits' data.
Dry weight of the eye lens (milligrams) versus age (days).
Plot of the local polynomial smoothing (dashed blue), the convex estimate (solid red), and data points (magenta).
The smoothing step is based on local polynomials of degree 1 (left) and degree 2 (right).
The fit obtained is really excellent, with no essential difference between degree 1 and 2.
}
\label{F:rabbits}
\end{figure}

\subsection{Two dimensional simulations}

In this section we briefly illustrate the finite sample properties of the convex estimate of a regression function in two dimensions by means of a simulation study.
For this purpose we considered the following convex regression function:
\[ f(x_1,x_2) = \max\left\{2 x_1^2 + x_2^2/2,3x_1+x_2\right\}, \]
which is convex, and only Lipschitz.
In \tref{Figure}{F:regression2d} we show the level curves of two estimated regression functions and the exact one in two simulations.
We took uniform grids of $10\times10$, and $20\times20$ in each situation for the explanatory variable, and added normal error with $\sigma=0.1$ to the value of $f(x_1,x_2)$ to emulate an observed variable.
The level curves shown in the figure show a very good fit, even for a coarse grid of only $20\times20$ points.

\begin{figure}
\capstart
\begin{center}
\includegraphics[width=.49\textwidth]{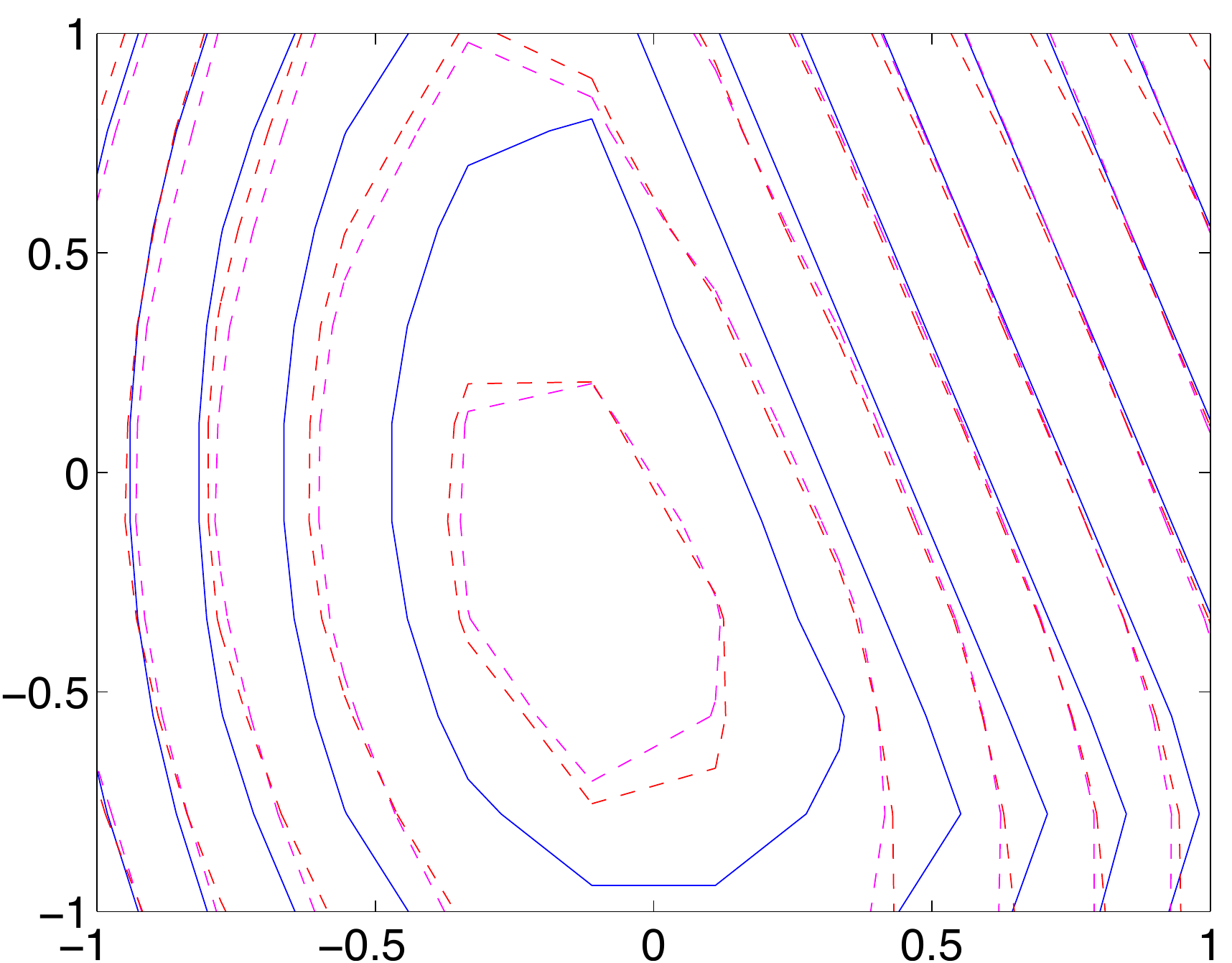}
\includegraphics[width=.49\textwidth]{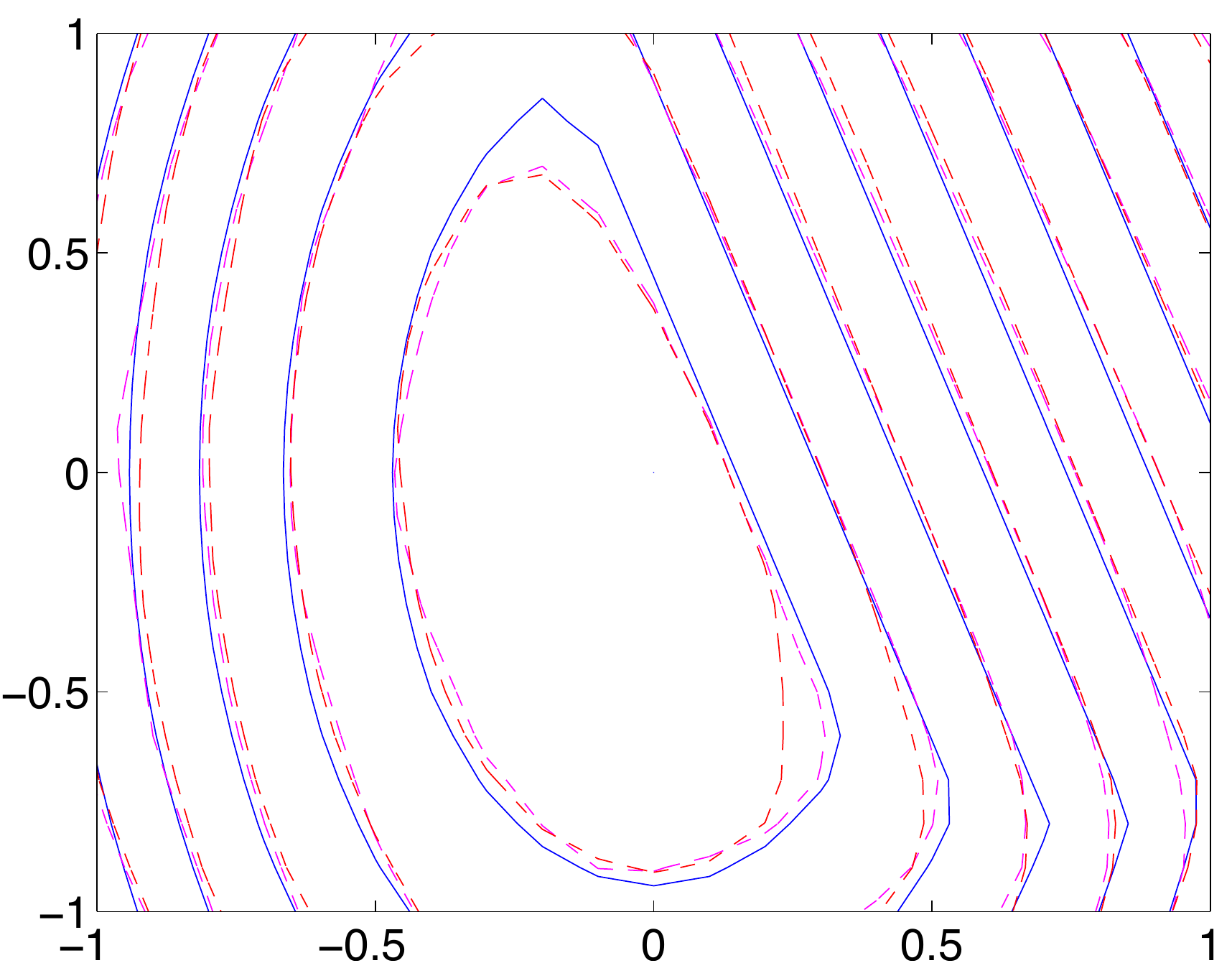}
\end{center}
\caption{\small
Level curves of two estimated regression functions (dashed red and magenta) and the exact one (solid blue) in two simulations with uniform design for the explanatory variables.
One for a grid of $10\times10$ (left) points and another for a grid of $20\times20$ (right).
The fit looks very good, even for a grid of only $20\times20$ points.}
\label{F:regression2d}
\end{figure}

In the second part of this simulation study we investigated the mean square error, bias and variance of our convex estimate.
For this we considered again the same two dimensional regression function $f$ and calculated by 2000 simulation runs the surfaces for the mean square error, squared bias and variance.
The results depicted in \tref{Figure}{F:varbiasmse2d} show that the variance is concentrated on the boundary but is one order of magnitude smaller than the squared bias and the mean square error.
These last two quantities are concentrated on the region of the domain where the regression function is not $C^1$.

\begin{figure}
\capstart
\begin{center}
Var \hfil\hfil Bias$^2$ \hfil\hfil MSE

\includegraphics[width=.32\textwidth]{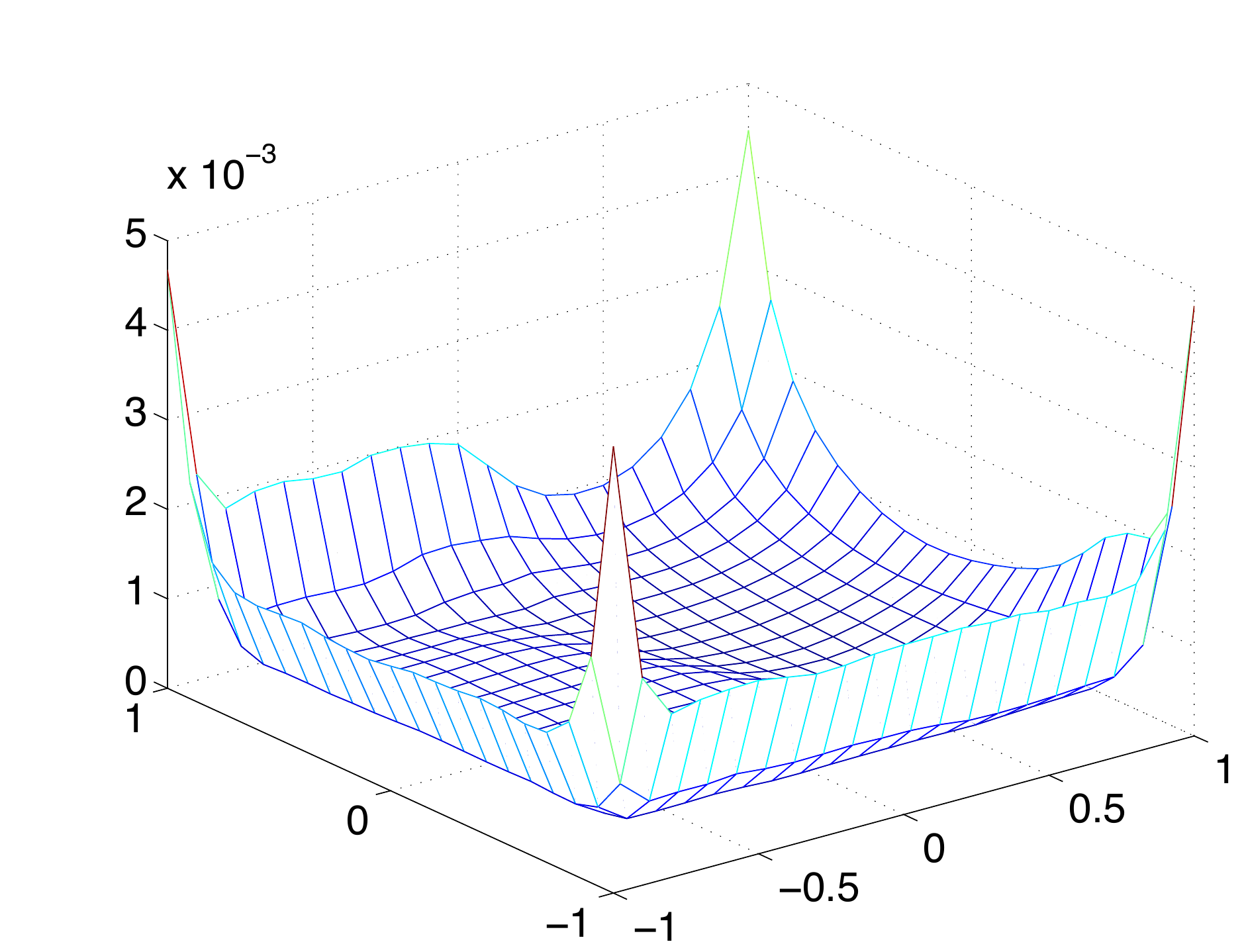}
\includegraphics[width=.32\textwidth]{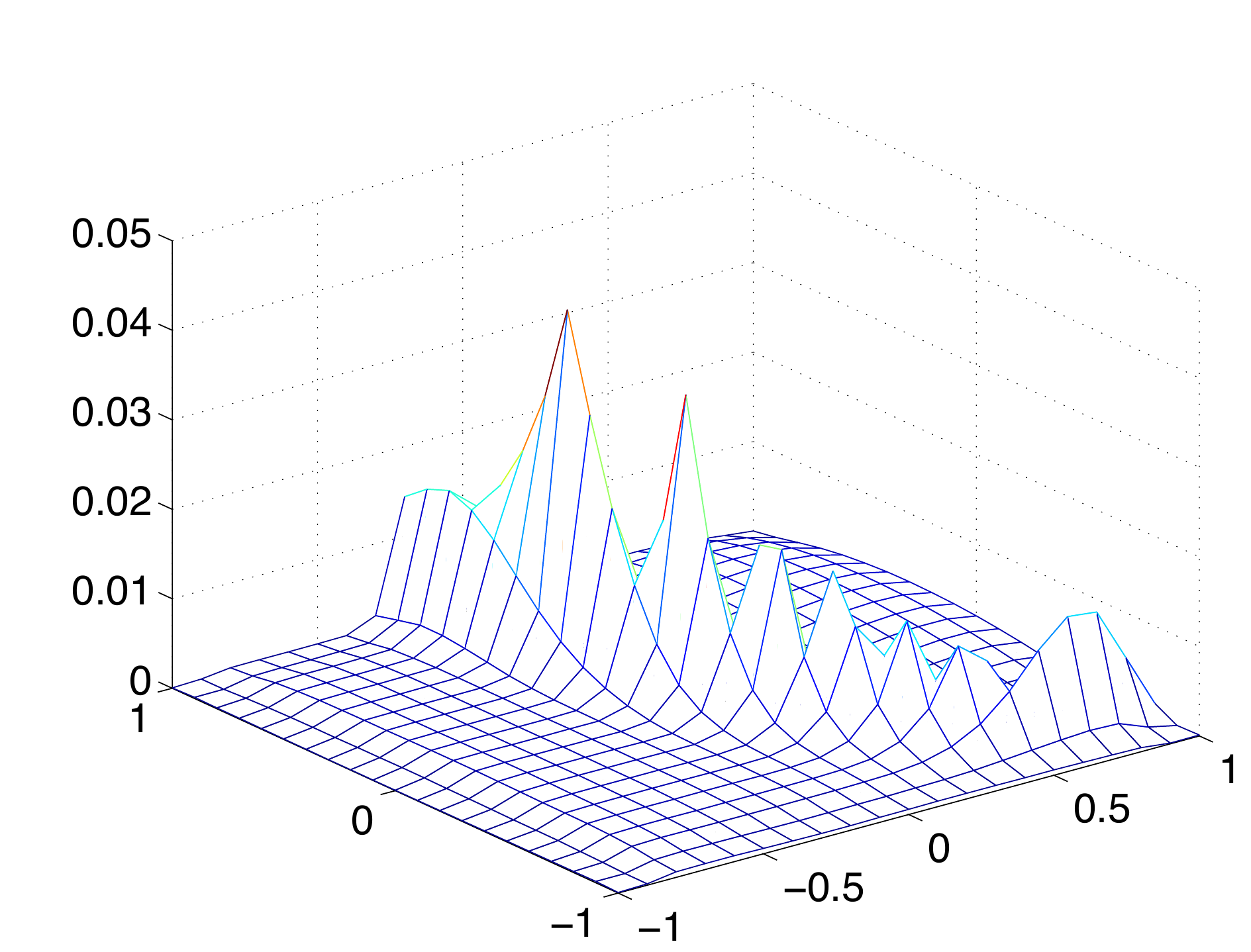}
\includegraphics[width=.32\textwidth]{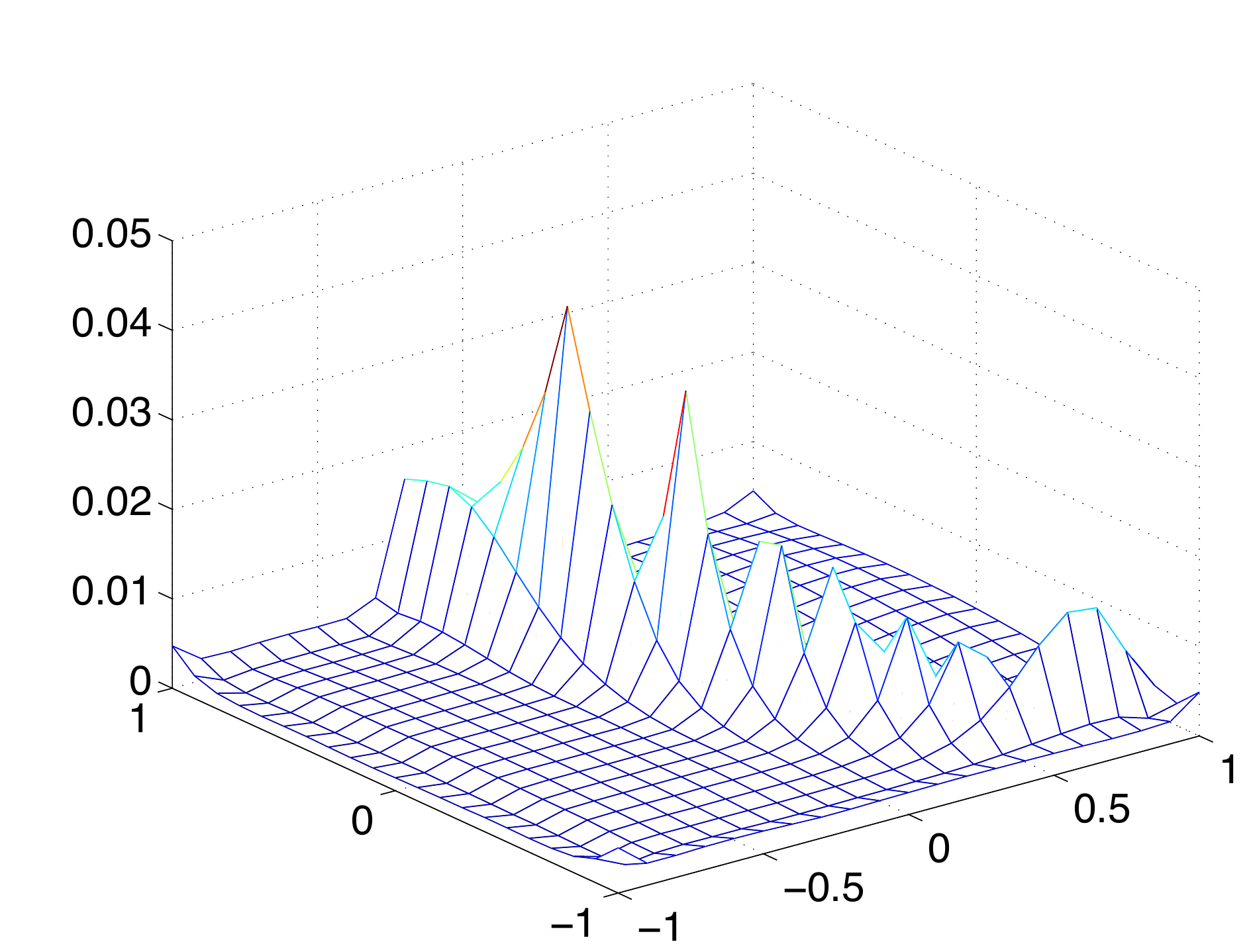}
\end{center}
\caption{\small
Variance (left), squared Bias (middle) and Mean Square Error (right) for the two dimensional simulation.
These indicators were obtained with 2000 simulation runs, using $10\times10$ (left) and $20\times20$ uniformly distributed design points for the explanatory variable, and normal error with $\sigma = 0.1$ for the observed variable.}
\label{F:varbiasmse2d}
\end{figure}

\subsection{Radiotherapy data}\label{sec:radio}

We studied a two dimensional example considered in~\citet{SHH05} (see also~\citealp{HSHKH05}),
who approximated the Pareto surface of a multiobjective optimization problem arising in the computation of the precise radiation dose.
This Pareto surface is convex under certain conditions, and it should be computed from some Pareto points that can be measured from the patient.
We obtained data from a patient of the Radboud University Nijmegen Medical Centre, in Nijmegen, the Netherlands.
The data correspond to a multiobjective optimization problem with three objectives and contains 69 data points, which, due to measuring errors, are not convex.

By using our method we are able to smooth the data,
obtaining a convex Pareto surface defined as a maximum of planes.
This surface is initially defined on the convex hull of the $\xx$ data, and we have extended it to a rectangular domain by considering the same maximum of planes.

In \tref{Figure}{F:example2d} we show the data points together with the convex regression surface (left), and the contours of the convex regression (right), showing an excelent fit of the data (see also Figure~2 in ~\citealp{HSHKH05}).

\begin{figure}
\capstart
\centering
\includegraphics[width=.49\textwidth]{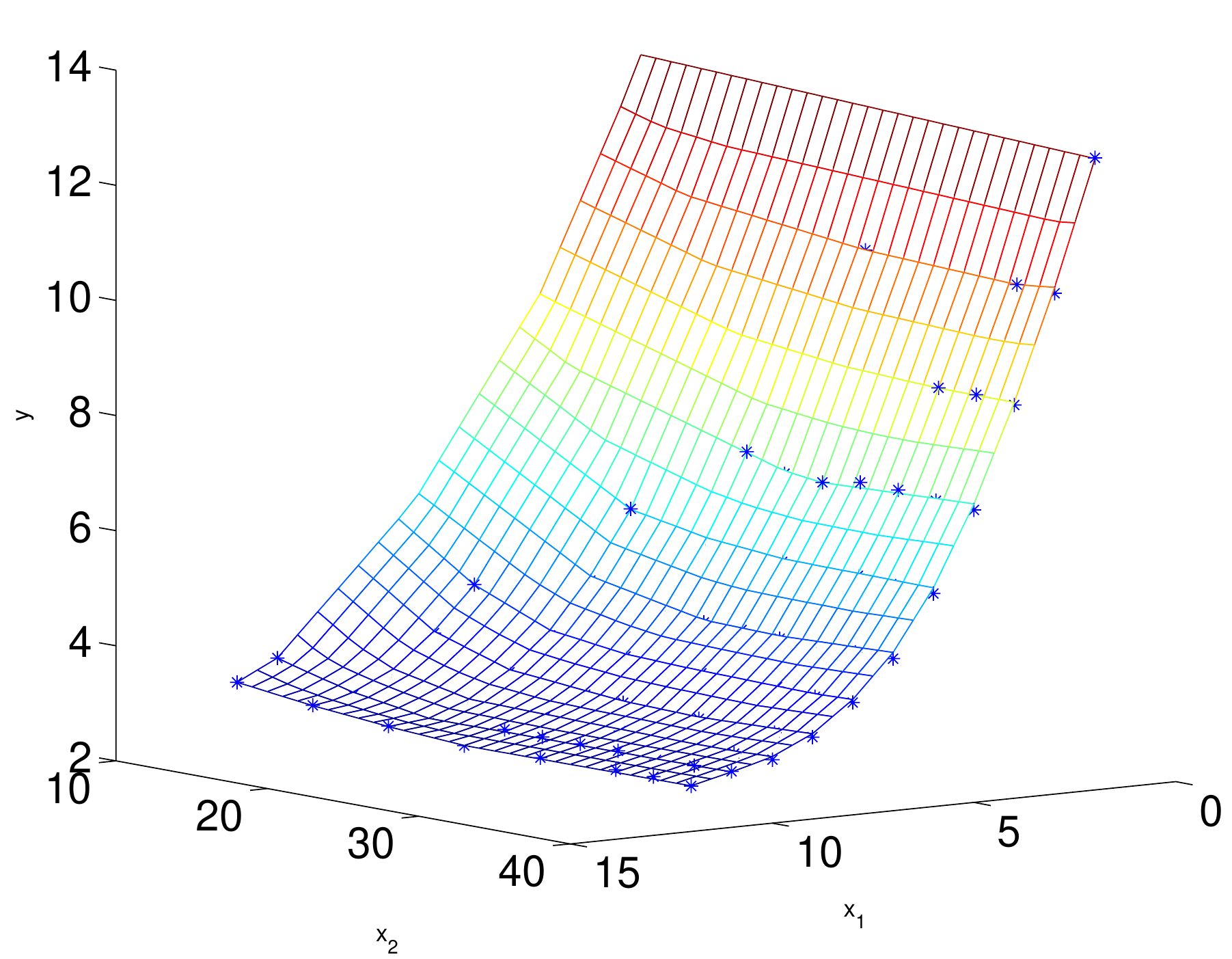}
\hfil
\includegraphics[width=.49\textwidth]{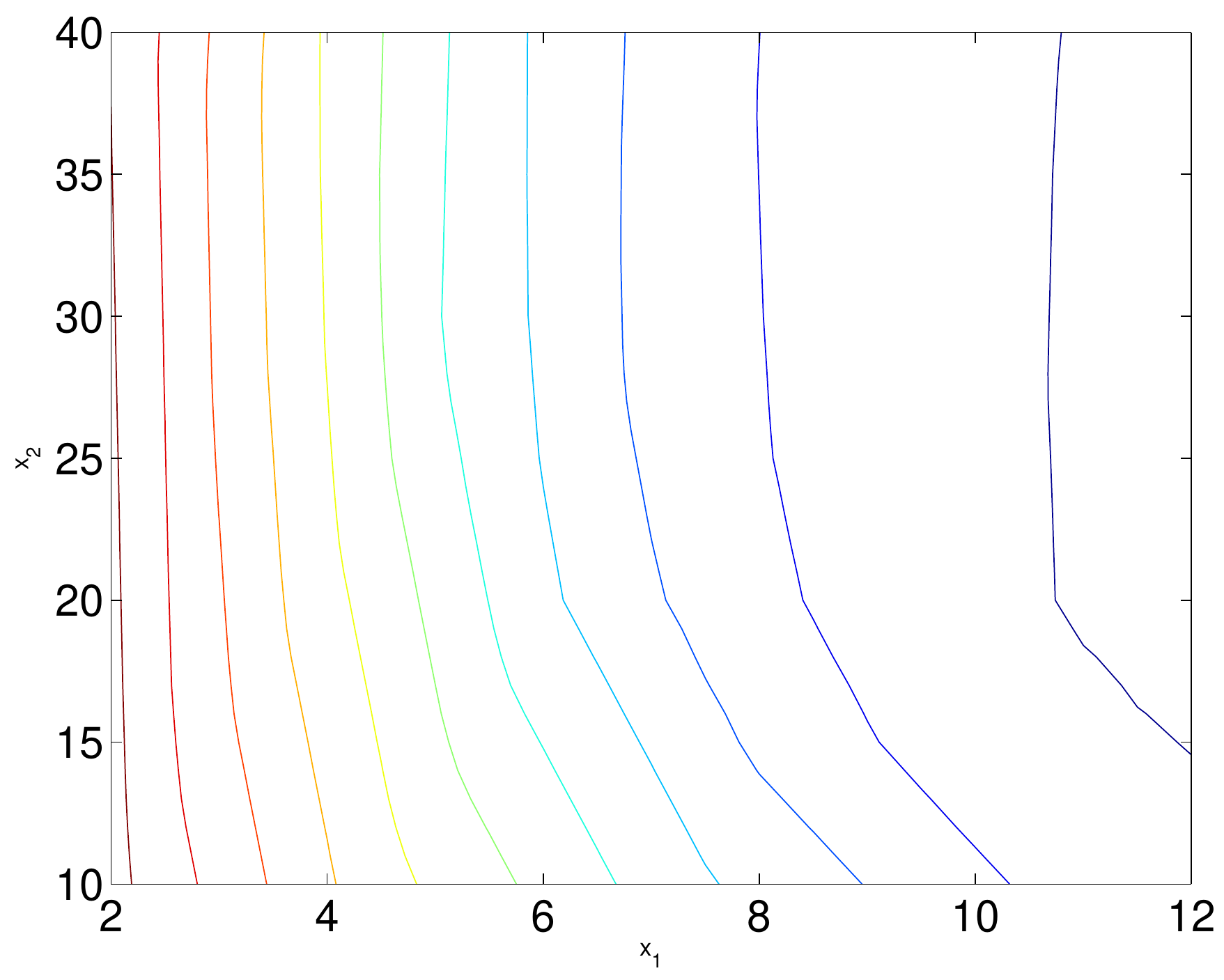}
\caption{Pareto surface obtained as the convex regression of the data points (left).
Contour curves of the convex graph (right) showing an excelent fit of the data (see also Figure~2 of ~\citealp{HSHKH05}).}
\label{F:example2d}
\end{figure}

\subsection*{Acknowledgements}

We would like to thank A.~Hoffmann and A.~Siem for sharing with us the data of the Radboud University Nijhmegen Medical Centre, used in \tref{Section}{sec:radio}.

\ifthenelse{\boolean{hyper}}{\phantomsection}{}
\protect\addcontentsline{toc}{section}{\refname}
\bibliography{convex}

{
\medskip
\noindent\rule{.5\textwidth}{1.5pt}\par
\medskip
\noindent
Corresponding author:\par
Liliana Forzani\par
Address: IMAL, Güemes 3450, 3000 Santa Fe, Argentina\par
e-mail: \url{liliana.forzani@gmail.com}
}

\end{document}